\documentclass{article}

\usepackage[english]{babel}
\usepackage{amsmath}
\usepackage{amsfonts}
\usepackage{graphicx}

\def\qed{{\hfill{\vrule height5pt width3pt depth0pt}\medskip}}

\newcommand{\cover}[1]{\stackrel{#1}{\Longrightarrow}}

\newcommand{\invcover}[1]{\stackrel{#1}{\Longleftarrow}}

\newcommand{\inte }{{\rm int}\,}

\newcommand{\bd }{\partial}

\def\comment#1{{}}

\def\makebold#1{\expandafter\def\csname bold#1\endcsname{{\mathbf #1}}}
\def\defboldone#1{\ifx#1\end \let\next=\relax
    \else \makebold#1 \let\next=\defboldone
    \fi \next}
\defboldone abcdefghijklmnopqrstuvwxyzABCDEFGHIJKLMNOPQRSTUVWXYZ\end

\newcommand{\boldlambda}{\boldsymbol\lambda}
\newcommand{\boldepsilon}{\boldsymbol\varepsilon}

\newcommand{\RectDef}[2]{[-#1,#1]\times[-#2,#2]}
\newcommand{\Rect}[2]{N_{#1,#2}}
\newtheorem{thm}{Theorem}[section]
\newtheorem{cor}[thm]{Corollary}
\newtheorem{lem}[thm]{Lemma}
\newtheorem{con}[thm]{Conjecture}
\newtheorem{defn}[thm]{Definition}
\newtheorem{rem}[thm]{Remark}

\newcommand{\dom }{\mbox{dom} \,}
\languageshorthands{english}

\begin{document}

\begin{center}
{\bf \large   Heteroclinic Connections between Periodic Orbits in
Planar Restricted Circular Three Body Problem - A Computer
Assisted Proof}

\vskip 0.5cm
{\large Daniel Wilczak} \\
  WSB -- NLU, Faculty of Computer Science,\\
  Department of Computational Mathematics,\\
  ul. Zielona 27, 33-300 Nowy S\c{a}cz, Poland\\
  e-mail: dwilczak@wsb-nlu.edu.pl

\vskip 0.2cm

and

\vskip 0.2cm
 {\large Piotr Zgliczynski}\footnote{Research supported in part by Polish KBN grant 2 P03A 011 18, 
         NSF grant DGE-98-04459  
         and PRODYN
         } \\
Jagiellonian University, Institute of Mathematics, \\ Reymonta 4,
30-059 Krak\'ow, Poland
\\ e-mail: zgliczyn@im.uj.edu.pl

\vskip 0.5cm

\today

\vskip 0.5cm

\end{center}

\begin{abstract}
The restricted circular three-body problem is considered for the
following parameter values $C=3.03$, $\mu=0.0009537$ - the values
for {\em Oterma} comet in the Sun-Jupiter system.
 We present a computer assisted proof of an existence of homo- and heteroclinic
 cycle between two Lyapunov orbits and an existence of symbolic
 dynamics on four symbols built on this cycle.
\end{abstract}

\section{Introduction and statement of results}

In paper \cite{KLMR} methods of dynamical system theory were used
(see also \cite{LR}) to explain rapid transitions from
heliocentric orbits  outside the orbit of Jupiter to heliocentric
orbits inside the orbit of Jupiter and vice versa for Jupiter
comets {\em Oterma} and {\em Gehrels 3}. To model this problem
authors in \cite{KLMR} used planar circular restricted three-body
problem and established that for a parameters corresponding to
Sun-Jupiter-Oterma system rapid transitions of Oterma are
explained by transversal intersections of stable and unstable
manifolds of two periodic orbits around libration points $L_1$ and
$L_2$. In fact the existence of symbolic dynamics on three symbols
was claimed.

The goal of this paper is develop and test tools which allow with
computer assistance to prove the results claimed in \cite{KLMR}.

Before we state our main results we  give a short description of
the planar restricted circular three-body problem. We follow the
paper \cite{KLMR} and  use the notation introduced there. Let $S$
and $J$ be two bodies called Sun and Jupiter, of masses
$m_s=1-\mu$ and $m_j=\mu$, $\mu \in (0,1)$, respectively. They
rotate in the plane in circles counter clockwise about their
common center and with angular velocity normalized as one. Choose
a rotating coordinate system (synodical coordinates) so that
origin is at the center of mass and the Sun and Jupiter are fixed
on the $x$-axis at $(-\mu,0)$ and $(1-\mu,0)$ respectively. In
this coordinate frame the equations of motion of a massless
particle called the comet or the spacecraft under the
gravitational action of Sun and Jupiter are (see \cite{KLMR} and
references given there)
\begin{equation}
  \ddot{x} - 2\dot{y}=\Omega_x(x,y), \qquad
  \ddot{y} + 2\dot{x}=\Omega_y(x,y),  \label{eq:PCR3BP}
\end{equation}
where
\begin{eqnarray*}
  \Omega(x,y)=\frac{x^2 + y^2}{2} + \frac{1 - \mu}{r_1} + \frac{\mu}{r_2}
     + \frac{\mu(1-\mu)}{2} \\
   r_1=\sqrt{(x+\mu)^2 + y^2}  , \qquad r_2=\sqrt{(x-1+\mu)^2 + y^2}
\end{eqnarray*}
Equations (\ref{eq:PCR3BP}) are called the equations of the planar
circular restricted three-body problem (PCR3BP). They have a first
integral called the {\em Jacobi integral}, which is given by
\begin{equation}
  C(x,y,\dot{x},\dot{y})= - (\dot{x}^2 + \dot{y}^2) +
  2\Omega(x,y).
\end{equation}

We  consider PCR3BP on the hypersurface
\begin{equation}
  \mathcal{M}(\mu,C)=\{ (x,y,\dot{x},\dot{y}) \: | \: C(x,y,\dot{x},\dot{y})=C
  \},
\end{equation}
and we restrict our attention to the following parameter values
$C=3.03$, $\mu=0.0009537$ - the parameter values for {\em Oterma}
comet in the Sun-Jupiter system (see \cite{KLMR}).

\begin{figure}
 \label{fig:HillsRegion}
 \centerline{\includegraphics{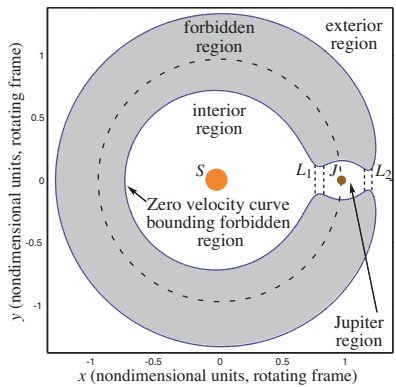}}
 \caption{Hills region for  PCR3BP with $C=3.03$, $\mu=0.0009537$ from  \cite{KLMR}.}
\end{figure}

The projection of $\mathcal{M}(\mu,C)$ onto position space is
called a Hill's region and gives the region in the $(x,y)$-plane,
where the comet is free to move. The Hill's region for the
parameter considered in this paper is shown on Figure
\ref{fig:HillsRegion} in white, the forbidden region is shaded.
The Hill's region consists of three regions: an interior (Sun)
region, an exterior region and Jupiter region.

In \cite{KLMR} a very good numerical evidence was given for the
following facts for the Sun-Jupiter-Oterma system:
\begin{description}
\item[0.] an existence of Lyapunov orbits $L_1^*$ and $L_2^*$ around
libration points $L_1$ and $L_2$, respectively. Both orbits are
hyperbolic and  are located in the Jupiter region.
\item[1.]  There exists a transversal heteroclinic orbit
connecting $L_1^*$ and $L_2^*$. There exists a transversal
     heteroclinic orbit connecting $L_2^*$ and $L_1^*$. Both
     orbits are in the Jupiter region. These orbits were discovered
     for the first time in \cite{KLMR}.
\item[2.] there exists a transversal homoclinic orbit to $L_1^*$
   in the interior (Sun) region,
\item[3.] there exists a transversal homoclinic orbit to $L_2^*$
   in the exterior region.
\end{description}
By transversal hetero- and homoclinic orbit, we mean that
appropriate unstable and stable manifolds intersect transversally.
For example in assertion {\bf 1}: the stable manifold of $L_1^*$
intersect transversally the unstable manifold of $L_2^*$.

It is now standard in dynamical system theory (see \cite{KLMR} and
references given there) to derive from assertions [0-3] an
existence of symbolic dynamics on four symbols $S, L_1^*, L_2^*,
X$ which the following allowed transition
\begin{displaymath}
  S \to S,L_1^*,  \quad L_1^* \to L_1^*,S, L_2^* \quad L_2^* \to
  L_1^*, L_2^*, X, \quad X \to X, L_2^*.
\end{displaymath}
In \cite{KLMR}( section 1.4 ) an existence of symbolic dynamics on
three symbols, only, was claimed. Instead of two symbols $L_1^*$
and $L_2^*$ one symbol $J$ for Jupiter region was used.

From the point of view of rapid transition of Oterma from interior
region to an exterior region and vice versa, an existence of
heteroclinic orbits between $L_1^*$ and $L_2^*$ claimed in
assertion 1 was of a special  importance, as they are an
indication of an existence of a dynamical channel joining an
interior with an exterior region.

  The following two theorems summarize main results of our paper
\begin{thm}
\label{thm:exists_homo_heteroclinic}
 For  PCR3BP with $C=3.03$, $\mu=0.0009537$  there exist two
 periodic solutions in the Jupiter region, $L_1^*$ and $L_2^*$ , called Lyapunov orbits,
 and there exists  heteroclinic connections
 between them, in both directions. Moreover for both orbits
 $L^*_1$ and $L_2^*$ there exists a homoclinic orbit in interior
 and exterior region, respectively.
\end{thm}

The next theorem  says that the homo- and heteroclinic connections
whose existence is established in Theorem
\ref{thm:exists_homo_heteroclinic} are topologically transversal
i.e. give rise to the symbolic dynamics, just as in the case of an
existence of transversal intersections of stable and unstable
manifolds.
\begin{thm}
\label{thm:symb_dyn_1}
 For  PCR3BP with $C=3.03$, $\mu=0.0009537$  there exist a
 symbolic dynamics on four symbols $\{S, X, L_1, L_2 \}$ corresponding to
 Sun and exterior regions and vicinity of  $L_1$ and $L_2$,
 respectively.
\end{thm}
A precise statement of this theorem with all necessary details
about the symbolic dynamics is given as Theorem
\ref{thm:symb_dyn_2} in Section \ref{sec:symbdyn}.

Hence we have proved assertion {\bf 0}, but we did'nt proved
assertions {\bf 1-3}, as we did'nt checked that stable and stable
manifolds of Lyapunov orbits intersect transversally. Instead we
had  proved that there is enough topological transversality
present to build a symbolic dynamics on it. The use of topological
tools was essential for the success of this work, as the rigorous
computation of stable and unstable manifolds appears to be much
more difficult than the computations reported in this paper.

Contents of our paper may be described as follows. In
Section~\ref{sec:pcr3bp} we continue our brief description of
PCR3BP, we define suitable Poincar\'e maps  and state their
symmetry properties. In Section~\ref{sec:topmeth} we present the
main topological tool used in this paper - the notion of covering
relation. In Section~\ref{sec:hyper} we describe how to link a
local hyperbolic behavior with covering relations to obtain homo-
and heteroclinic orbits. In Sections~\ref{sec:Lyap_orb} and
\ref{sec:homo_hetero} we report the results of our rigorous
computations for PCR3BP and  we prove
Theorem~\ref{thm:exists_homo_heteroclinic}. In
Section~\ref{sec:symbdyn} we show how to use symmetries of PCR3BP
together with covering relations to complete the proof of
Theorem~\ref{thm:symb_dyn_1}. Section~\ref{sec:num} contains the
details of the  numerical part of proof, we mainly discuss the
question of an efficient approach to a verification of covering
relations.  We also include all initial data, so that a willing
reader with his own code can verify our claims. In
Section~\ref{sec:programm} we provide a brief description of basic
classes used in our programm and we give the names of functions,
which perform proofs of numerical lemmas, in order to facilitate
an understanding of the source code, which we had made available
on-line (see\cite{W}). In Section~\ref{sec:future} we discuss some
natural extensions of our results.

What is new in this paper besides giving a proof for some results
from \cite{KLMR}? First of all it shows how to successfully link
numerically cheap $C^0$-methods (the covering relations) with much
more numerically expensive $C^1$-methods (a local hyperbolicity).
This was previously done for maps (see \cite{GZ}), only. The main
obstacle in applications to ODEs was the lack of an efficient
$C^1$ ODE solver. Such a solver - called a $C^1$-Lohner algorithm
- was recently proposed by the second author in \cite{Z3}.

Other novelties in this paper are some theoretical improvements in
the theory of covering relations. We  use  an abstract definition
of covering relation from \cite{GiZ} and we show how to use a
symmetry which involves a time reflection (for Poincar\'e map this
correspond to taking an inverse map) together with covering
relations. Both these improvements  (in numerical algorithms and
in theory of covering relations) allow to reduce the computation
time considerably (probably by two orders of magnitude). The total
computation time on a PC using a 1.1Ghz Celeron processor is less
than 40 minutes.

\section{Properties of PCR3BP: Poincar\'e maps and symmetries}
\label{sec:pcr3bp}
 In this section we continue our brief description of PCR3BP
we started in the Introduction and we introduce various notations
which will be used throughout the paper.

The PCR3BP has three unstable collinear equilibrium points on the
Sun-Jupiter line, called $L_1$, $L_2$ and $L_3$ (see Fig. 1.4 in
\cite{KLMR} ), whose eigenvalues include one real and one
imaginary pair. The value of $C$ (Jacobi constant) at the point
$L_i$ will be denoted by $C_i$. An linearization at $L_i$ for
$i=1,2$ for the parameter range considered here, shows that these
points are of center-saddle type (see \cite{KLMR}). By theorem of
Moser \cite{M} it follows that for $C < C_i$ and $|C-C_i|$ small
enough, there exist hyperbolic periodic orbits, $L_i^*$, around
$L_i$, called Lyapunov orbits. Observe that for a fixed value of
$C < C_i$ an existence of Lyapunov orbit $L_i^*$ is not settled by
the Moser theorem and has to be proved.

As was mentioned in the Introduction we restrict our attention to
the following parameter values $C=3.03$, $\mu=0.0009537$ - the
parameter values for {\em Oterma} comet in the Sun-Jupiter system
(see \cite{KLMR}). Since we work with fixed parameter values we
usually drop the dependence of various objects defined throughout
the paper on $\mu$ and $C$, so for example
$\mathcal{M}=\mathcal{M}(\mu,C)$. For our parameter values we have
$C_2 > C > C_3$ ( this means that we are considering Case 3 from
section 3.1 in \cite{KLMR}).

 We consider  Poincar\'e
sections: $\Theta = \{(x,y,\dot{x},\dot{y})\in \mathcal{M} \:| \:
y=0\}$, $\Theta_+ = \Theta\cap\{\dot{y}>0\}$, $\Theta_- =
\Theta\cap\{\dot{y}<0\}$.

On $\Theta_\pm$ we can express $\dot{y}$ in terms of $x$ and
$\dot{x}$ as follows
\begin{equation}
\dot{y} = \pm\sqrt{2\Omega(x,0) - \dot{x}^2-C}
\end{equation}
Hence the sections $\Theta_\pm$ can be parameterized by two
coordinates $(x,\dot{x})$ and we will use this identification
throughout the paper. More formally, we have the transformation
$T_{\pm}:\mathbb{R}^2 \to \Theta_{\pm}$
 given by the following formula
\begin{equation}
 T_\pm(x,\dot{x})=(x,0,\dot{x},\pm\sqrt{2\Omega(x,0) - \dot{x}^2-C}\mbox{ })
\end{equation}
The domain of $T_\pm$ is given by an inequality $2\Omega(x,0) -
\dot{x}^2-C \geq 0$.

Let $\pi_{\dot{x}}:\Theta_\pm \longrightarrow\mathbb{R}$ and
$\pi_{x}:\Theta_\pm \longrightarrow\mathbb{R}$ denote the
projection onto $\dot{x}$ and $x$ coordinate, respectively. We
have $\pi_{\dot{x}}(x_0,\dot{x}_0)=\dot{x}_0$ and
$\pi_{x}(x_0,\dot{x}_0)=x_0$.

We will say that $(x,\dot{x}) \in \Theta_\pm$ meaning that
$(x,\dot{x})$ represents two-dimensional coordinates of a point on
$\Theta_\pm$. Analogously we  give a meaning to the statement $M
\subset \Theta_\pm$ for a set $M \subset \mathbb{R}^2$.

We define the following Poincar\'e maps between sections
\begin{eqnarray*}
  P_+: \Theta_+ \to \Theta_+ \\
  P_-: \Theta_- \to \Theta_- \\
  P_{\frac{1}{2},+}: \Theta_+ \to \Theta_- \\
  P_{\frac{1}{2},-}: \Theta_- \to \Theta_+.
\end{eqnarray*}
As a rule the sign $+$ or $-$ tells that the domain of the maps
$P_\pm$ or $P_{\frac{1}{2},\pm}$ is contained in $\Theta_\pm$ (the
same sign). Observe that
\begin{displaymath}
   P_+(x)=P_{\frac{1}{2},-} \circ P_{\frac{1}{2},+}(x), \qquad
   P_-(x)=P_{\frac{1}{2},+} \circ P_{\frac{1}{2},-}(x)
\end{displaymath}
whenever $P_+(x)$ and $P_-(x)$ are defined. These identities
express the following simple fact: to return to $\Theta_+$ we need
to cross $\Theta$ with negative $\dot{y}$ (this is
$P_{\frac{1}{2},+}$ first and then we return to $\Theta$ with
$\dot{y} >0$ (this is $P_{\frac{1}{2},-}$).

Sometimes we will drop signs in $P_\pm$ and $P_{\frac{1}{2},\pm}$,
hence $P(z)=P_{+}(z)$ if $z \in \Theta_+$ and $P(z)=P_{-}(z)$ if
$z \in \Theta_-$, a similar convention will be applied to
$P_{\frac{1}{2}}$.

\subsection{Symmetry properties of PCR3BP}
\label{subsec:symPCR3BP}
 Notice that PCR3BP has the following symmetry
\begin{equation}
R(x,y,\dot{x},\dot{y},t)=(x,-y,-\dot{x},\dot{y},-t),
\end{equation}
 which
expresses the following fact, if $(x(t),y(t))$ is a trajectory for
PCR3BP, then $(x(-t),-y(-t))$ is also a trajectory for PCR3BP.
From this it follows immediately that
\begin{eqnarray}
  \mbox{if} \quad P_\pm(x_0,\dot{x}_0)=(x_1,\dot{x}_1), \qquad \mbox{then} \quad
    P_\pm(x_1,-\dot{x}_1)=(x_0,-\dot{x}_0) \label{eq:sym_Pf} \\
  \mbox{if} \quad P_{\frac{1}{2},\pm}(x_0,\dot{x}_0)=(x_1,\dot{x_1}), \qquad \mbox{then} \quad
    P_{\frac{1}{2},\mp}(x_1,-\dot{x}_1)=(x_0,-\dot{x}_0) \nonumber
\end{eqnarray}

We will denote also by $R$ the map $R:\Theta_\pm \to \Theta_\pm$
$R(x,\dot{x})=(x,-\dot{x})$ for $(x,\dot{x}) \in \Theta_\pm$. Now
eq. (\ref{eq:sym_Pf}) can be written as
\begin{eqnarray}
 \mbox{if} \quad P_\pm(x_0)=x_1, \qquad \mbox{then} \quad
    P_\pm(R(x_1))=R(x_0) \label{eq:sym_PfS} \\
  \mbox{if} \quad P_{\frac{1}{2},\pm}(x_0)=x_1, \qquad \mbox{then} \quad
    P_{\frac{1}{2},\mp}(R(x_1))=R(x_0) \nonumber
\end{eqnarray}

\section{Topological tools}
\label{sec:topmeth}

In this section we present  main topological tools used in this
paper. The crucial notion is that of {\em a covering relation}.
This notion in various forms was introduced  in papers
\cite{Z0,Z1,Z2,Z4}. Here we follow the most recent and most
general version introduced in \cite{GiZ} and the reader is
referred there for proofs.

\subsection{h-sets}

{\bf Notation:} For a given norm in $\mathbb{R}^n$ by $B_n(c,r)$
we will denote an open ball of radius $r$ centered at $c \in
\mathbb{R}^n$. When the dimension $n$ is obvious from the context
we will drop the subscript $n$. Let $S^n(c,r)=\partial
B_{n+1}(c,r)$, by the symbol $S^n$ we will denote $S^n(0,1)$. We
set $\mathbb{R}^0=\{0\}$, $B_0(0,r)=\{0\}$, $\partial
B_0(0,r)=\emptyset$.

For a given set $Z$, by $\inte Z$, $\overline{Z}$, $\partial Z$ we
denote the interior, the closure and the boundary of $Z$,
respectively. For the map $h:[0,1]\times Z \to \mathbb{R}^n$ we
set $h_t=h(t,\cdot)$.  By $\mbox{Id}$ we denote an identity map.
For a map $f$, by $\mbox{dom} (f)$ we will denote the domain of
$f$. Let $f : \Omega \subset {\mathbb R}^n \to {\mathbb R}^n$ a
continuous map we will say that $X \subset \dom(f^{-1})$ if the
map $f^{-1}:X \to {\mathbb R}^n$ is well defined and continuous.

\begin{defn}
\label{def:covrel} A $h$-set, $N$, is the  object consisting of
the following data
\begin{itemize}
 \item $|N|$ - a compact subset of ${\mathbb R}^n$, a support of
 $N$
 \item $u(N),s(N) \in \{0,1,2,\dots\}$, such that $u(N)+s(N)=n$
 \item a homeomorphism $c_N:{\mathbb R}^n \to
   {\mathbb R}^n={\mathbb R}^{u(N)} \times {\mathbb R}^{s(N)}$,
     such that
      \begin{displaymath}
        c_N(|N|)=\overline{B_{u(N)}}(0,1) \times
        \overline{B_{s(N)}}(0,1).
      \end{displaymath}
\end{itemize}
We set
\begin{eqnarray*}
  N_c=\overline{B_{u(N)}}(0,1) \times \overline{B_{s(N)}}(0,1), \\
   N_c^-=\partial \overline{ B_{u(N)}}(0,1) \times
\overline{B_{s(N)}}(0,1) \\
N_c^+=\overline{B_{u(N)}}(0,1) \times
\partial \overline{B_{s(N)}}(0,1) \\
  N^-=c_N^{-1}(N_c^-) , \quad N^+=c_N^{-1}(N_c^+)
\end{eqnarray*}
\end{defn}

Hence a $h$-set, $N$, is a product of two closed balls in some
coordinate system. The numbers, $u(N)$ and $s(N)$, stand for the
dimensions of nominally unstable and stable directions,
respectively. The subscript $c$ refers to the new coordinates
given by homeomorphism $c_N$. We will call $N^-$ ($N_c^-$) an exit
set of $N$ and $N^+$ ($N^+_c$) an entry set of $N$. Observe that
if $u(N)=0$, then $N^-=\emptyset$ and if $s(N)=0$, then
$N^+=\emptyset$.

\begin{defn}
Let   $N$ be a $h$-set. We define a $h$-set $N^T$ as follows
\begin{itemize}
 \item $|N^T|=|N|$
 \item $u(N^T)=s(N)$,  $s(N^T)=u(N)$
 \item We
 define a homeomorphism $c_{N^T}:{\mathbb R}^n \to   {\mathbb R}^n={\mathbb R}^{u(N^T)} \times {\mathbb
R}^{s(N^T)}$,
 by
      \begin{displaymath}
        c_{N^T}(x)= j(c_{N}(x)) ,
      \end{displaymath}
      where $j: {\mathbb R}^{u(N)} \times {\mathbb R}^{s(N)} \to {\mathbb R}^{s(N)} \times {\mathbb R}^{u(N)}$
      is given by $j(p,q)=(q,p)$.
\end{itemize}
\qed
\end{defn}
Observe that $N^{T,+}=N^-$ and $N^{T,-}=N^+$. This operation is
useful in the context of inverse maps, as it was first pointed out
in \cite{A}.

\subsection{Covering relations}

For $n>0$ and a continuous map $f:S^n \to S^n$ by $d(f)$ we denote
the degree of $f$ \cite{DG}. For $n=0$ we define the degree,
$d(f)$,  as follows. Observe first that $S^0=\{-1,1\}$. We set
\begin{equation}
  d(f)=
  \begin{cases}
    1, & \text{if $f(1)=1$ and $f(-1)=-1$}, \\
    -1, & \text{if $f(1)=-1$ and $f(-1)=1$}, \\
    0,  & \text{otherwise}.
  \end{cases}
\end{equation}

\begin{defn}
Assume $n > 0$. Let $f:\overline{B_n}(0,1) \to  {\mathbb R}^n$,
such that
\begin{equation}
  0 \notin  f(\partial B(0,1))).
\end{equation}
We define a map $s_f:S^{n-1} \to S^{n-1}$ by
\begin{equation}
  s_f(x) = \frac{f(x)}{\|f(x) \|}.
\end{equation}
\end{defn}

\begin{defn}
Assume $N,M$ are $h$-sets, such that $u(N)=u(M)=u$ and
$s(N)=s(M)=s$. Let $f:|N| \to {\mathbb R}^n$ be continuous. Let
$f_c= c_M \circ f \circ c_N^{-1}: N_c \to {\mathbb R}^u \times
{\mathbb R}^s$. Let $w$ be a nonzero integer. We say that
\begin{displaymath}
  N\cover{f,w} M
\end{displaymath}
($N$ $f$-covers $M$ with degree $w$) iff the following conditions
are satisfied
\begin{description}
\item[1.] there exists a continuous homotopy $h:[0,1]\times N_c \to {\mathbb R}^u \times {\mathbb R}^s$,
   such that the following conditions hold
   \begin{eqnarray}
      h_0&=&f_c,  \label{eq:hc1} \\
      h([0,1],N_c^-) \cap M_c &=& \emptyset  \label{eq:hc2} \\
      h([0,1],N_c) \cap M_c^+ &=& \emptyset \label{eq:hc3}
   \end{eqnarray}
\item[2.1] If $u > 0$, then there exists a  map $A:{\mathbb R}^u \to {\mathbb R}^u$, such that
   \begin{eqnarray}
    h_1(p,q)&=&(A(p),0), \quad \mbox{where $p \in {\mathbb R}^u$ and $q \in {\mathbb
    R}^s$}\\
      A(\partial B_u(0,1)) & \subset & {\mathbb R}^u \setminus
      \overline{B_u(0,1)}  \label{eq:mapaway}
   \end{eqnarray}
  Moreover, we require that
\begin{displaymath}
  d(s_A)=w,
\end{displaymath}
\item[2.2] If $u=0$, then
\begin{eqnarray}
  h_1(x) &=& 0  \qquad \mbox{for $x \in N_c$} \\
  w &=& 1.
\end{eqnarray}
\end{description}
\end{defn}

Intuitively,  $N \cover{f} M$ if $f$ stretches  $N$ in the
'nominally unstable' direction, so that its projection onto
'unstable' direction in $M$ covers in topologically nontrivial
manner projection of $M$. In the 'nominally stable' direction $N$
is contracted by $f$. As a result $N$ is mapped across $M$ in the
unstable direction, without touching $M^+$. An example of covering
relation on the plane with one unstable direction is shown on
Figure~\ref{fig:cov}.

\begin{defn}
 Assume $N,M$ are $h$-sets, such that $u(N)=u(M)=u$ and
$s(N)=s(M)=s$.  Let $g:\Omega \subset {\mathbb R}^n \to {\mathbb
R}^n$. Assume that $g^{-1}:|M| \to {\mathbb R}^n$ is well defined
and continuous. We say that $N \invcover{g,w} M$ ($N$
$g$-backcovers $M$ with degree $w$) iff $M^T \cover{g^{-1},w}
N^T$.
\end{defn}

The following theorem proved in \cite{GiZ} is one of  main tools
used in this paper. Various versions of this theorem (without
backcovering) using slightly weaker notions of covering relations
or even without  an explicitly defined notion of covering relation
were given in \cite{Z0,Z1,Z2,Z4}. In the planar case this theorem
with backcovering  was stated also in \cite{A}.
\begin{thm}
\label{th:top}
 Assume $N_i$, $i=0,\dots,k$, $N_k=N_0$ are
$h$-sets and for each $i=1,\dots,k$ we have either
\begin{equation}
  N_{i-1} \cover{f_i,w_i} N_{i} \label{eq:dirgcov}
\end{equation}
or  $|N_{i}| \subset \dom(f_i^{-1})$ and
\begin{equation}
  N_{i-1} \invcover{f_i,w_i} N_{i}.  \label{eq:invgcov}
\end{equation}

Then there exists a point $x \in \inte |N_0|$, such that
\begin{eqnarray}
   f_i \circ f_{i-1}\circ \cdots \circ f_1(x) &\in& \inte |N_i|, \quad i=1,\dots,k \\
  f_k \circ f_{k-1}\circ \cdots \circ f_1(x) &=& x
\end{eqnarray}
\end{thm}

Obviously we cannot make any claim about the uniqueness of $x$ in
Theorem~\ref{th:top}.

\subsection{Covering relation on the plane with one nominally expanding direction ($u=1$)}
In this section we discuss the case, when $u=s=1$, hence we have
only one nominally expanding and one nominally contracting
direction. The basic idea here is: the set $N^-$ consists from two
disjoint components and all possible values of the degree $w$ in
covering relation are $\pm1$. This allows to give sufficient
conditions for an existence of covering relations, which are
relatively easy to verify.

\begin{defn}
Let $N$ be a h-set, such that $u(N)=s(N)=1$. We set
\begin{eqnarray*}
  N_c^{le}&=&\{-1\} \times [-1,1] \\
  N_c^{re}&=&\{1\} \times [-1,1] \\
  S(N)_c^{l} &=& (-\infty,-1) \times {\mathbb R} \\
  S(N)_c^{r} &=& (1,\infty) \times {\mathbb R}.
\end{eqnarray*}
We define
\begin{eqnarray*}
  N^{le}=c_N^{-1} (N_c^{le}), \quad  N^{re}=c_N^{-1} (N_c^{re}),
  \\
  \quad  S(N)^{l}=c_N^{-1} (S(N)^{l}), \quad  S(N)^{r}=c_N^{-1}
  (S(N)^{r}).
\end{eqnarray*}
We will call $N^{le}$, $N^{re}$, $S(N)^l$ and $S(N)^r$ the left
edge, the right edge, the left side and right side of $N$,
respectively.
\end{defn}

It is easy to see that $N^-=N^{le} \cup N^{re}$.

The triple $(|N|,\overline{S(N)^l}, \overline{S(N)^r})$ is a t-set
from \cite{AZ}. As in \cite{AZ} we will use the following notation
for $S(N)^{r,l}$.

\begin{displaymath}
  N^l=\overline{S(N)^l}, \qquad  N^r={\overline{S(N)^r}}
\end{displaymath}

\begin{figure}
[ptb]
\begin{center}
\includegraphics[
height=1.51226in, width=2.61726in
]%
{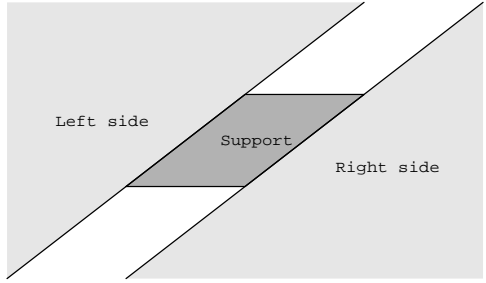}%
\caption{An example of h-set on the plane.}%
\label{pic:magicset}%
\end{center}
\end{figure}

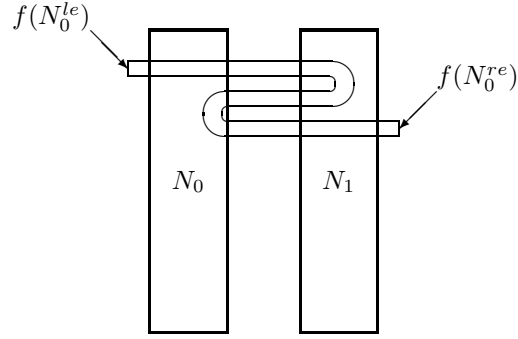
\begin{figure}[ptb]
\setlength{\unitlength}{1cm}
\begin{picture}(11,6)(-4,-1)
\thicklines \put(0,0){\framebox(1,4){$N_0$}}
\put(2,0){\framebox(1,4){$N_1$}} \thinlines
\put(-0.8,4){\vector(1,-1){0.5}}
\put(-0.8,4){\makebox(0,0)[br]{$f(N_0^{le})$}}
\put(-0.3,3.6){\line(1,0){2.5}} \put(-0.3,3.4){\line(1,0){2.5}}
\put(-0.3,3.4){\line(0,1){0.2}} \put(2.2,3.3){\oval(1,0.6)[r]}
\put(2.2,3.3){\oval(0.5,0.2)[r]} \put(2.2,3.2){\line(-1,0){1}}
\put(2.2,3.0){\line(-1,0){1}} \put(1.2,2.9){\oval(1,0.6)[l]}
\put(1.2,2.9){\oval(0.5,0.2)[l]} \put(1.2,2.6){\line(1,0){2.1}}
\put(1.2,2.8){\line(1,0){2.1}} \put(3.3,2.6){\line(0,1){0.2}}
\put(3.8,3.2){\vector(-1,-1){0.5}}
\put(3.8,3.2){\makebox(0,0)[bl]{$f(N_0^{re})$}}
\end{picture}
\caption{An example of an $f-$covering relation: $N_{0}
\Longrightarrow N_{0},N_{1}$}
 \label{fig:cov}
\end{figure}

\begin{rem}
For all h-sets used in this paper  the support is a parallelogram.
A usual picture of a h-set is given in Figure \ref{pic:magicset}.

A typical picture illustrating covering relation on the plane with
one 'unstable' direction is given on Figure~\ref{fig:cov}.
\end{rem}

The following theorem was proved in \cite{GiZ} for any $n>1$ and
$u(N)=1$. Here we rewrite it for the planar case in a slightly
different notation (we use $N^l$ and $N^r$ for $S(N)^l$ and
$S(M)^r$, respectively).
\begin{thm}
\label{thm:suffcovu1} Let $n=2$ and let $N$, $M$ be two h-sets in
${\mathbb R}^n$, such that $u(N)=u(M)=1$ and $s(N)=s(M)=1$. Let
$f:|N| \to {\mathbb R}^n$ be continuous.

Assume that there exists $q_0 \in {\overline B}_s(0,1)$, such that
following conditions are satisfied
   \begin{eqnarray}
      f(c_N([-1,1] \times \{q_0 \})) &\subset& \inte (M^l  \cup |M| \cup M^r )
               \label{eq:imfu1} \\
      f(|N|) \cap M^+ &=& \emptyset,  \label{eq:imfu3}
   \end{eqnarray}
and one of the following two conditions holds
\begin{eqnarray}
  f(N^{le}) \subset M^{l} \quad \mbox{and} \quad f(N^{re}) \subset
  M^{r} \label{eq:covorpr} \\
  f(N^{le}) \subset M^{r} \quad \mbox{and} \quad f(N^{re}) \subset
  M^{l} \label{eq:covorrev}.
\end{eqnarray}

 Then there exists  $w = \pm 1$, such that
\begin{displaymath}
  N\cover{f,w} M
\end{displaymath}
\end{thm}

\subsection{Representation of the h-sets}
\label{subsec:reptsets} In this paper we  use very simple h-sets,
namely the support is a parallelogram.

A h-set is  defined by specifying the triple $N=t(c,u,s)$, where
$c,u,s \in \mathbb{R}^2$, such that $u,s$ are linearly
independent. We set
\begin{eqnarray*}
  |N|& =& \{ x \in \mathbb{R}^2 \ | \ \exists_{t_1,t_2 \in [-1,1]} \quad  x=c + t_1 u + t_2 s
  \} \\
   &=& c + [-1,1] \cdot u + [-1,1] \cdot s \\
   c_N(t_1,t_2)&=& c + t_1 u + t_2 s
\end{eqnarray*}
 In this representation $c$ is a center point of
the parallelogram $N$, $u$ represents an oriented half-length in
the 'unstable' direction and $s$ is an oriented half-length in the
'stable' direction. See Fig.~\ref{pic:magicset} for an example of
h-set in this representation.

 We have
\begin{eqnarray*}
  N^{le} &=& c - u + [-1,1] \cdot s  \\
  N^{re} &=& c + u + [-1,1] \cdot s \\
  N^{l} &=& c + (-\infty,-1] u + (-\infty,\infty) s \\
  N^{r} &=& c + [1,\infty) u + (-\infty,\infty) s.
\end{eqnarray*}
We introduce  notions of  top and bottom edges of $N$, $N^{te}$
and $N^{be}$  by
\begin{eqnarray*}
  N^{be} &=& c + [-1,1]\cdot u -  s  \\
  N^{te} &=& c + [-1,1]\cdot u +  s
\end{eqnarray*}

Let us recall that the symmetry $R:\mathbb{R}^2 \to \mathbb{R}^2$
introduced in section~\ref{subsec:symPCR3BP} was given by
\begin{displaymath}
  R(x_1,x_2)=(x_1,-x_2)
\end{displaymath}

\begin{defn}
\label{def:symtset} A h-set, $N$, will be called {\em an
$R$-symmetric h-set} if $N=t(c,u,s)$ for some $c,u,s \in
\mathbb{R}^2$, such that
\begin{eqnarray*}
  R(c)&=&c \\
  R(u)=s \quad &\mbox{or}& \quad R(u)=-s
\end{eqnarray*}
\end{defn}

\begin{figure}[ptb]
\begin{center}
\includegraphics[width=2.5in]{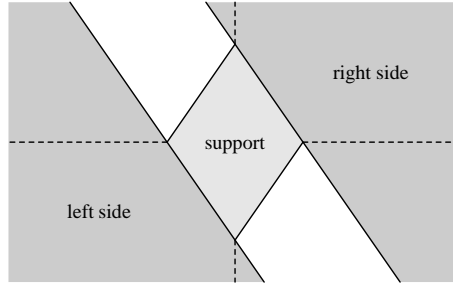}
\caption{An example of an $R$-symmetric h-set} \label{pic:symset}
\end{center}
\end{figure}

Figure~\ref{pic:symset} shows an example of a $R$-symmetric h-set.
Symmetry properties of such h-set are apparent.

\subsection{Action of $R$ on h-sets}
The symmetry of $P_{1/2,\pm}$ and $P_\pm$ expressed in
(\ref{eq:sym_Pf}) relates the maps and theirs inverse, hence
beside mapping the support of $N$ by $R$ it will switch also the
nominally stable and unstable directions. This motivates the
following definition of the action of the symmetry $R$ on h-sets
\begin{defn}
\label{def:rhsets}
  Let $N$ be a h-set. We define a h-set $R(N)$ as follows
  \begin{itemize}
    \item $|R(N)|=R(|N|)$
    \item $u(R(N))=s(N)$ and $s(R(N))=u(N)$
    \item the homeomorphism $c_{R(N)}:\mathbb{R}^n \to \mathbb{R}^{u(R(N))} \times \mathbb{R}^{u(S(N))}$
        is given by
\begin{displaymath}
   c_{R(N)}= c_{N^T} \circ R^{-1}
\end{displaymath}
  \end{itemize}
  \qed
\end{defn}

Observe that according to above definition we have
\begin{eqnarray}
R(N)^-=R(N^+)=R(N^{T,-}) \nonumber \\
R(N)^+=R(N^-)=R(N^{T,+}) \nonumber \\
R(t(c,u,s))= t(R(c),R(s),R(u)) \label{eq:rnahset}
\end{eqnarray}

As an immediate consequence of equation (\ref{eq:rnahset}) we
obtain
\begin{lem}
\label{lem:sym_on_tset} Let $N=t(c,u,s)$ be an $R$-symmetric
h-set. Then  $R(N)=N$.
\end{lem}

We have the following easy lemma
\begin{lem}
\label{lem:covsym} Let $f_i:\dom(f_i) \subset \mathbb{R}^2 \to
 \mathbb{R}^2$ for $i=1,2$  continuous and invertible on some open sets. Assume that
\begin{equation}
  \mbox{if} \quad f_1(x)=x_1, \qquad \mbox{then} \qquad
  f_2(Rx_1)=R(x).
\end{equation}
Let $N_1$, $N_2$ be h-sets such that
\begin{displaymath}
  N_1 \cover{f_1} N_2.
\end{displaymath}
Then
\begin{displaymath}
  R(N_2) \invcover{f_2} R(N_1)
\end{displaymath}
\end{lem}
{\em Proof:} From Def.~\ref{def:rhsets} and the assumed symmetry
it follows immediately that
\begin{equation}
   R(N_1)^T \cover{f_2^{-1}} R(N_2^T).
\end{equation}
 \qed

From above Lemma and (\ref{eq:sym_PfS}) we obtain
\begin{cor}
\label{cor:covsym} If $N_1 \cover{P_\pm} N_2$, then $R(N_2)
\invcover{P_\pm} R(N_1)$.

If $N_1 \cover{P_{\frac{1}{2},+}} N_2$, then $R(N_2)
\invcover{P_{\frac{1}{2},-}} R(N_1)$.

If $N_1 \cover{P_{\frac{1}{2},-}} N_2$, then $R(N_2)
\invcover{P_{\frac{1}{2},+}} R(N_1)$.
\end{cor}

\section{$C^1$ tools}
\label{sec:hyper}

The goal of this section is to describe the tools which allow in
the presence of hyperbolic fixed points for a map  to prove an
existence of homo- and heteroclinic trajectories.

In this section we recall  the results from \cite{GZ} with some
additions (see also \cite{WZ} where the method was outlined for
the first time). In the symbol of covering relation we will drop
the degree part, hence we will use $N \cover{f} M$ instead of $N
\cover{f,w} M$ for some nonzero $w$.

\subsection{General theorems}
Let $P:{\mathbb R}^n \to {\mathbb R}^n$ be a $C^1$-map.  For any
set $X$ we define an interval  matrix  $[DP(X)] \subset {\mathbb
R}^{n\times n}$ to be an interval enclosure of $DP(X)$ given by
\begin{displaymath}
  M \in [DP(X)] \quad \mbox{iff} \quad  \inf_{x \in X}DP(x)_{ij}
  \leq  M_{ij} \leq \sup_{x \in X}DP(x)_{ij} \quad i,j=1,2,\dots,n
\end{displaymath}

\begin{lem}
\label{lem:error} Let $N$ be a convex set. Assume $x_0, x_1 \in
N$. Then
\begin{equation}
   P(x_1) - P(x_0)  \in [DP(N)] \cdot (x_1 - x_0).
   \label{eq:basicerror}
\end{equation}
Moreover,  there exists a matrix $M \in [DP(N)]$ such that
\begin{equation}
  P(x_1) - P(x_0)  = M \cdot (x_1 - x_0)
\end{equation}
\end{lem}
{\em Proof:}
\begin{eqnarray*}
P(x_1) - P(x_0)= \int_0^1 \frac{d P}{d t} (x_0 + t(x_1 - x_0)) dt
= \\
  \int_0^1 \frac{\partial P}{\partial x} (x_0 + t(x_1 - x_0)) dt \cdot (x_1 - x_0)
\end{eqnarray*}
To finish the proof observe that
$$
M=\int_0^1 \frac{\partial P}{\partial x} (x_0 + t(x_1 - x_0)) dt
\in [DP(N)]
$$
\qed

Let $Id:\mathbb{R}^n \to \mathbb{R}^n$ denote the identity map.
\begin{thm}
\label{thm:ewunique} Let $N$ be a convex set. Assume that
\begin{displaymath}
  0 \notin \mbox{det}([DP(N)]-Id) = \{ t=det(M-Id) \ | \ M \in [DP(N)]
  \},
\end{displaymath}
then $N$ contains at most one fixed point of $P$.
\end{thm}
{\em Proof:} Assume that $x_0,x_1 \in N$ are fixed points for $P$.
Then from lemma \ref{lem:error} it follows that
\begin{equation}
  x_1 - x_0 = P(x_1) - P(x_0) = M \cdot (x_1 - x_0)
\end{equation}
for some matrix $M \in [DP(N)]$. Hence $x_1 - x_0$ is in the
kernel of $M - Id$.  From our assumption it follows that
$\mbox{det}(M - Id) \neq 0$, hence $x_1=x_0$. \qed

Consider a two--dimensional function $f(x)=(f_1(x),f_2(x))^T$,
where $x=(x_1,x_2)^T$. We assume that $f(0)=0$, i.e. $0$ is a
fixed point of $f$. For  a convex set $U$, such that $0 \in U$ we
define intervals $\boldlambda_1(U)$, $\boldepsilon_1(U)$,
$\boldepsilon_2(U)$ and  $\boldlambda_2(U)$ by
\begin{equation}
Df(U)=\left(
\begin{array}{cc}
\boldlambda_1(U) & \boldepsilon_1(U) \\
\boldepsilon_2(U) & \boldlambda_2(U)
\end{array}\right).
\end{equation}

Since $f(0)=0$ then from Lemma \ref{lem:error} it follows that
\begin{eqnarray*}
f_1(x)\in\boldlambda_1(U) x_1 + \boldepsilon_1(U) x_2 \\
f_2(x)\in\boldepsilon_2(U) x_1 + \boldlambda_2(U) x_2
\end{eqnarray*}

Let
\begin{eqnarray*}
\varepsilon_1'(U)=\sup\{|\varepsilon|:\varepsilon\in\boldepsilon_1(U)\},
\quad
& \varepsilon_2'(U)=\sup\{|\varepsilon|:\varepsilon\in\boldepsilon_2(U)\}, \\
\lambda_1'(U)=\inf\{|\lambda_1|:\lambda_1\in\boldlambda_1 (U)\},
\quad &
\lambda_2'(U)=\sup\{|\lambda_2|:\lambda_2\in\boldlambda_2(U) \}.
\end{eqnarray*}

Let us define the rectangle $\Rect{\alpha_1}{\alpha_2}$ by
\begin{displaymath}
\Rect{\alpha_1}{\alpha_2}=\RectDef{\alpha_1}{\alpha_2}.
\end{displaymath}

\begin{defn}\cite[Def. 1]{GZ}
\label{def:hyperbolic}
 Let $x_*$ be a fixed point for the map $f$.
We say that $f$ is {\em hyperbolic} on $N \ni x_*$, if there
exists a local coordinate system on $N$, such that in this
coordinate system
\begin{eqnarray}
x_*&=&0 \\
 \label{eq:hyper} \varepsilon'_1(N) \varepsilon'_2(N) &<& (1 -
\lambda'_2(N))(\lambda'_1(N) - 1).  \\
N&=&\Rect{\alpha_1}{\alpha_2},
\end{eqnarray}
 where $\alpha_1 >0$, $\alpha_2 > 0$ are such
that the following conditions are satisfied
\begin{eqnarray}
\frac{\varepsilon'_1(N)}{\lambda_1'(N) -1} <
\frac{\alpha_1}{\alpha_2} < \frac{1 -
\lambda'_2(N)}{\epsilon'_2(N)}. \label{eq:hratio}
\end{eqnarray}
\end{defn}

It is easy to see that that for the map $f$ to be hyperbolic on $N$ it is necessary that
$\lambda'_1 > 1$,$\lambda'_2 < 1$ and the linearization of $f$ at $x_*$ is hyperbolic with one
stable and unstable direction.

\begin{thm}\cite[Thm. 3]{GZ}
\label{th:hyper} Assume that $f$ is hyperbolic on $N$. Then
\begin{description}
\item[1.]  if $f^{k}(x) \in N$ for $k \geq 0$, then
 $ \lim_{k\to \infty}f^k(x)=x_*$,
\item[2.] if $y_k \in N$ and $f(y_{k-1})=y_k$
  for $k \leq 0$, then $\lim_{k\to -\infty}y_k=x_*$.
\end{description}
\end{thm}

The next theorem  shows how we can combine ${\mathcal C}^0$- and
${\mathcal C}^1$-tools to prove the existence of asymptotic orbits
with prescribed itinerary.
\begin{thm}\cite[Thm. 4]{GZ}
\label{th:covhomo} Assume that $g$ is hyperbolic on $N_m$ and $f$
hyperbolic on $N_0$. Let $x_g \in N_m$ be a fixed point for $g$
and $x_f \in N_0$ be a fixed point for $f$.
\begin{description}
\item[1.]
If
\begin{equation} \label{eq:seqcov1}
N_0 \cover{f_0} N_1 \cover{f_1} N_2 \cover{f_2} \dots
\cover{f_{m-1}} N_m \cover{g} N_m
\end{equation}
then there exists $x_0 \in N_0$ such that
\begin{eqnarray*}
f_{i-1} \circ f_{i-2} \circ \dots \circ f_0(x_0) \in N_{i}   \quad & \mbox{ for }\quad i=1,\dots,m, \\
g^k \circ f_{m-1} \circ  \dots \circ f_0(x_0) \in N_m        \quad & \mbox{ for }\quad k>0, \\
\lim_{k\to \infty} g^k \circ f_{m-1} \circ  \dots \circ
f_0(x_0)=x_g.
\end{eqnarray*}
\item[2.]
If
\begin{equation} \label{eq:seqcov2}
N_0 \cover{f} N_0 \cover{f_0} N_1 \cover{f_1} N_2 \cover{f_2}
\dots \cover{f_{m-1}} N_m
\end{equation}
then there exists a sequence $(x_k)_{k=-\infty}^{0}$,
$f(x_k)=x_{k+1}$  for $k<0$ such that
\begin{eqnarray*}
x_k \in N_0                                                \quad & \mbox{ for }\quad k\leq 0, \\
f_{i-1} \circ f_{i-2} \circ \dots \circ f_0(x_0) \in N_{i} \quad & \mbox{ for }\quad i=1,\dots,m, \\
\lim_{k\to -\infty} x_k=x_f.
\end{eqnarray*}
\item[3.]
If
\begin{equation} \label{eq:seqcov3}
N_0 \cover{f} N_0 \cover{f_0} N_1 \cover{f_1} N_2 \cover{f_2}
\dots \cover{f_{m-1}} N_m \cover{g} N_m
\end{equation}
then there exists a sequence $(x_k)_{k=-\infty}^{0}$,
$f(x_k)=x_{k+1}$ for $k<0$ such that
\begin{eqnarray*}
x_k \in N_0,  \quad k \leq 0, \\
f_{i-1} \circ f_{i-2} \circ \dots \circ f_0(x_0) \in N_{i} \quad & \mbox{ for }\quad i=1,\dots,m,\\
g^n \circ f_{m-1} \circ  \dots \circ f_0(x_0) \in N_m      \quad & \mbox{ for } \quad n>0,\\
\lim_{k\to -\infty} x_k=x_f,\\
\lim_{k\to \infty} g^k \circ f_{m-1} \circ  \dots \circ
f_0(x_0)=x_g.
\end{eqnarray*}
\end{description}
\end{thm}
The above theorem can be used without any modifications for
proving the existence of trajectories converging to periodic
orbits. In this case we  consider higher iterates of maps $f$ and
$g$ in (\ref{eq:seqcov1}), (\ref{eq:seqcov2}) and
(\ref{eq:seqcov3}).

\subsection{How to prove an existence of an heteroclinic orbit, fuzzy sets.}
\label{subsec:howtoprovehomo}

To prove an existence of an heteroclinic orbit we
 want to use the third assertion in Theorem \ref{th:covhomo} for
$g=f$, but in order to make the exposition easier to follow  we
use  two different maps $f$ and $g$. Observe that to apply this
theorem directly one needs to know an exact location of two fixed
points $x_f \in N_0$ and $x_g \in N_m$, because the sets $N_0$ and
$N_m$ are centered on $x_f$ and $x_g$ respectively. But exact
coordinates of $x_f$ and $x_g$ are usually unknown. We overcome
this obstacle in three steps as follows

\noindent {\bf 1. Finding very good estimates for $x_f$ and
$x_g$}.
 In this paper we use an argument based on symmetry to obtain tight bounds
 for $x_f$ and $x_g$. In \cite{GZ} a rigorous interval Newton
 algorithm was used. Let us denote by $M_f$ and $M_g$ obtained estimates for
 $x_f$ and $x_g$, respectively.

 We choose one fixed point  $x_f \in M_f$  and $x_g \in M_g$
 for further considerations.

\noindent {\bf 2. $C^1$-computations, hyperbolicity}  We choose a
set $U_f$, $M_f \subset U_f$, on which we  compute rigorously
$[Df(U_f)]$. Then we have to choose a coordinate system, in which
the matrix $[Df(U_f)]$ will be as close as possible to the
diagonal one. In this paper we have chosen  numerically obtained
stable and unstable eigenvectors. Let us denote these eigenvectors
by $u$ and $s$, where $u$ corresponds to unstable direction and
$s$ is pointing in the stable direction. Assume that this process
gives us a coordinate frame in which
\begin{equation}
 \varepsilon'_1(U_f) \varepsilon'_2(U_f) < (1 -
\lambda'_2(U_f))(\lambda'_1(U_f) - 1).  \label{eq:hyperU}
\end{equation}
From (\ref{eq:hyperU}) it follows easily that there exists
$\alpha_1 >0$, $\alpha_2>0$ such that
\begin{equation}
\frac{\varepsilon'_1(U_f)}{\lambda_1'(U_f) -1} <
\frac{\alpha_1}{\alpha_2} < \frac{1 -
\lambda'_2(U_f)}{\epsilon'_2(U_f)}, \label{eq:hratioU}.
\end{equation}
Observe that above inequality specifies only the ratio
$\alpha_1/\alpha_2$, hence we can find a pair
$(\alpha_1,\alpha_2)$ such that condition (\ref{eq:hratioU}) and
the following condition holds
\begin{equation}
  M_f + \alpha_1 \cdot [-1,1] \cdot u + \alpha_2 \cdot [-1,1] \cdot s \subset U_f
\end{equation}
We now define a h-set $N_0$ by
\begin{equation}
  N_0 = t(x_f,\alpha_1 u,\alpha_2 s).
\end{equation}
Obviously  $f$ is hyperbolic on $N_0$. Observe that the
hyperbolicity implies uniqueness of $x_f$ in $N_0$.

We do similar construction for $g$ to obtain $N_m=t(x_g,\beta_1
\bar{u},\beta_2 \bar{s})$.

\noindent {\bf 3. Covering relations for fuzzy h-sets}. We have to
verify the following covering relations
\begin{eqnarray}
  N_0 \cover{f} N_0 \cover{f} N_1 \label{eq:covh1}  \\
  N_{m-1} \cover{f_{m-1}} N_m \cover{g} N_{m} \label{eq:covh2}.
\end{eqnarray}
As was mentioned above we don't know the h-sets  $N_0$, $N_m$
explicitly, but we know that
\begin{eqnarray*}
  N_0 \in {\widetilde N}_0=\{ t(c,\alpha_1 u, \alpha_2 s) \ | \ c \in M_f   \} \\
  N_m \in {\widetilde N}_m=\{ t(c,\beta_1 \bar{u}, \beta_2 \bar{s}) \ | \ c \in M_g   \}.
\end{eqnarray*}
Above equations define {\em a fuzzy h-set}, as a collection of
h-sets. We can now extend the definition of covering relations to
fuzzy h-sets as follows.
\begin{defn}
\label{def:covfuzzy} Let $f$ be a continuous map on the plane.
Assume ${\widetilde N_1}$, ${\widetilde N_2}$ are fuzzy h-sets
(collections of h-set) and $R$ is a h-set.
\begin{itemize}
\item we say that ${\widetilde N_1} \cover{f} R$ iff  $M \cover{f}
R$ for all $M \in {\widetilde N_1}$.
\item we say that $R \cover{f} {\widetilde N_1}$ iff  $R \cover{f}
M$ for all $M \in {\widetilde N_1}$.
\item we say that ${\widetilde N_1} \cover{f} {\widetilde N_2}$ iff  $M_1 \cover{f}
M_2$ for all $M_1 \in {\widetilde N_1}$ and $M_2 \in {\widetilde
N_2}$.
\end{itemize}
\end{defn}

With the above definition is obvious that to prove the covering
relations in equations (\ref{eq:covh1}) and (\ref{eq:covh2}) it is
enough to show that
\begin{eqnarray}
 {\widetilde N}_0 \cover{f} {\widetilde N}_0 \cover{f} N_1   \\
  N_{m-1} \cover{f_{m-1}} {\widetilde N}_m \cover{g} {\widetilde N}_{m}.
\end{eqnarray}

In practice (in rigorous numerical computations) it is convenient
to think about a fuzzy h-set ${\widetilde N}$ as an parallelogram
with thickened edges, hence all tools developed to verify covering
relations for h-sets can be easily extended to fuzzy h-sets.

\section{The Lyapunov orbits.}
\label{sec:Lyap_orb}
 In this section we present a computer
assisted proof of an existence and hyperbolicity of the Lyapunov
orbits around libration points. Hence we realize here step 1 and 2
from section \ref{subsec:howtoprovehomo} on our way to the proof
of an existence of heteroclinic connection.

As in previous section in the symbol of covering relation we will
drop the degree part, hence we will use $N \cover{f} M$ instead of
$N \cover{f,w} M$ for some nonzero $w$.

\subsection{Existence of Lyapunov orbits.}

\begin{thm}
\label{thm:exists_per} Let $x_1=0.9208034913207400196$, $x_2=1.081929486841799903$.
\begin{itemize}
\item There exists a fixed point $L_1^*=(x_1^*,0) \in\Theta_+$ for $ P_+$, such that
\begin{equation}
    |x_1^*-x_1|<\eta_1=6\cdot10^{-14}\label{eq:inequal_1}
\end{equation}
\item There exists a fixed point $L_2^*=(x_2^*,0) \in\Theta_-$ for $P_-$ such that
\begin{equation}
    |x_2^*-x_2|<\eta_2=10^{-13}\label{eq:inequal_2}
\end{equation}
\end{itemize}
\end{thm}

\begin{figure}[hdtb]
 \centerline{\includegraphics{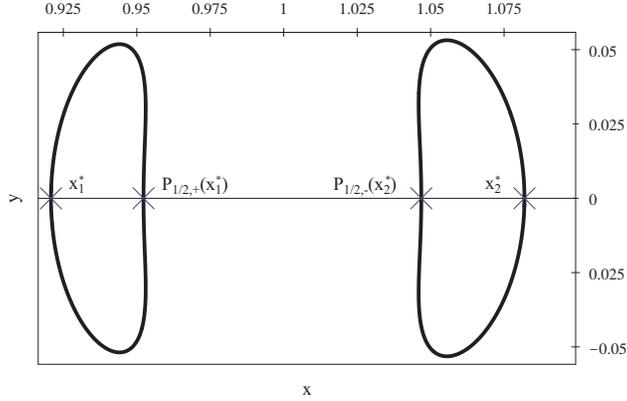}}
 \caption{The Lyapunov orbits and the location of  $x_i^*$.\label{fig:Lap_orb}
}
\end{figure}

{\em Proof:} We consider two intervals
$I_1:=[x_1-\eta_1,x_1+\eta_1] \times \{0\} \subset \Theta_+$,
$I_2:=[x_2-\eta_2,x_2+\eta_2] \times \{0\} \subset \Theta_-$. The
location of $x_i$ is schematically shown on Figure
\ref{fig:Lap_orb}.

Let us recall, that by
$\pi_{\dot{x}}:\Theta\longrightarrow\mathbb{R}$ we denote the
projection onto $\dot{x}$ coordinate. With a computer assistance
we proved the following
\begin{lem}
\label{lem:cov_rel_Lap} The maps
$P_{\frac{1}{2},+}:I_1\longrightarrow \Theta_-$ and
$P_{\frac{1}{2},-}:I_2\longrightarrow \Theta_+$ are well defined
and continuous. Moreover we have the following properties
\begin{eqnarray}
  \pi_{\dot{x}}(P_{\frac{1}{2},+}(x_1-\eta_1,0))<0,\qquad
  \pi_{\dot{x}}(P_{\frac{1}{2},+}(x_1+\eta_1,0))>0 \label{eq:Lapineq1} \\
  \pi_{\dot{x}}(P_{\frac{1}{2},-}(x_2-\eta_2,0))<0,\qquad
  \pi_{\dot{x}}(P_{\frac{1}{2},-}(x_2+\eta_2,0))>0 \label{eq:Lapineq2}
\end{eqnarray}
\end{lem}

\begin{figure}[h]
 \centerline{\includegraphics{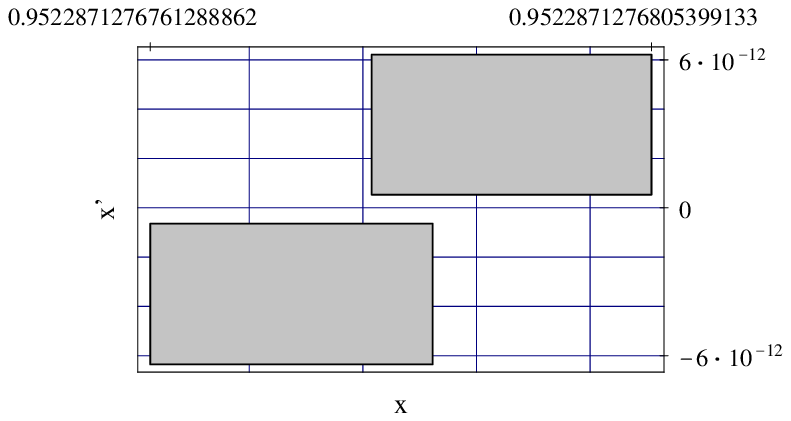}}
 \caption{Rigorous enclosure of $P_{\frac{1}{2},+}(x_1-\eta_1,0)$
 (a box in lower left corner ) and $P_{\frac{1}{2},+}(x_1+\eta_1,0)$ (a box in upper right corner)}
\label{fig:l1}
\end{figure}

\begin{figure}[h]
 \centerline{\includegraphics{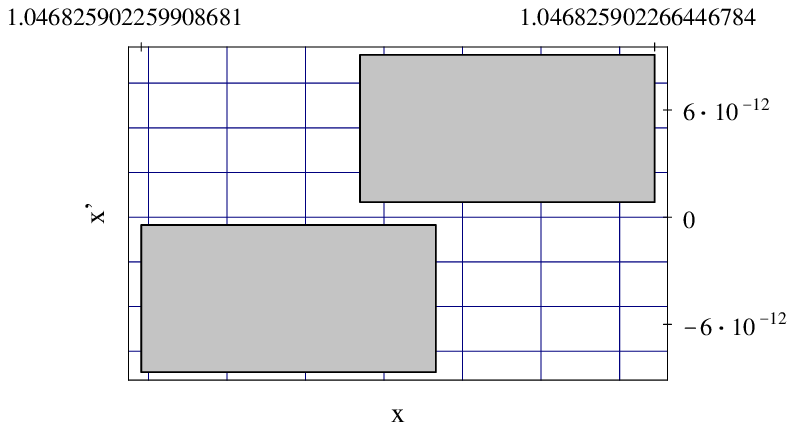}}
 \caption{Rigorous enclosure of $P_{\frac{1}{2},-}(x_2-\eta_2,0)$
 (a box in lower left corner ) and $P_{\frac{1}{2},-}(x_2+\eta_2,0)$ (a box in upper right corner )}
\label{fig:l2}
\end{figure}

Figures \ref{fig:l1} and \ref{fig:l2} display  rigorous enclosures
for $P_{\frac{1}{2},+}(x_1 \pm \eta_1,0)$ and
$P_{\frac{1}{2},-}(x_2 \pm \eta_2,0)$, respectively.

Now we are ready to finish the  proof of Theorem
\ref{thm:exists_per}. From Lemma \ref{lem:cov_rel_Lap} and the
Darboux property it follows that there exist points $x_1^*\in
\mbox{int}(I_1)$ and $x_2^*\in\mbox{int}(I_2)$, such that
\begin{eqnarray}
 P_{\frac{1}{2},+}(x_1^*,0)= (x_1^0,0)\\
 P_{\frac{1}{2},-}(x_2^*,0)= (x_2^0,0).
\end{eqnarray}
An application of symmetry properties of $P_{1/2,\pm}$ (see eq.
(\ref{eq:sym_Pf})) gives
\begin{eqnarray}
 P_+(x_1^*,0)= (x_1^*,0)\\
 P_-(x_2^*,0)= (x_2^*,0).
\end{eqnarray}
\qed

\subsection{Hyperbolicity in the neighborhood of Lyapunov orbits.}
\label{subsec:hyperLyap} The goal of this section is to prove that
$P$ is hyperbolic in the sense of Definition~\ref{def:hyperbolic}
in some neighborhood of points $L_1^*$ and $L_2^*$.

Let us define
\begin{eqnarray*}
    u_1  =  (1,2.5733011), \quad  s_1 = (-1,2.5733011),\\
    u_2  =  (1,2.2817915),\quad  s_2 = (-1,2.2817915).
\end{eqnarray*}
 These vectors appear to be a good
approximation for unstable ($u_i$)  and stable eigenvectors
($s_i$) at $L_i^*$ on the $(x,\dot{x})$-plane. Observe that
$R(u_i)=-s_i$, this is in agreement with symmetry of $P_\pm$
stated in equation (\ref{eq:sym_Pf}). We will also use $(u_i,s_i)$
later, as the coordinate directions for good coordinate frame in
the proof of hyperbolicity of $P_+$ and $P_-$ in the neighborhood
of $L^*_i$.

 Let  $H_i^1=t(h_i,u_i^1,s_i^1)$ and $H_i^2=t(h_i,u_i^2,s_i^2)$  for $i=1,2$
 denote  h-sets on the $(x,\dot{x})$ plane, where

\begin{eqnarray}
\begin{array}{ll}\label{eq:a}
  h_1 = (x_1,0),&
  h_2 = (x_2,0) \\
  \alpha_1=3 \cdot 10^{-10}, &
  \alpha_2=4\cdot10^{-10} \\
  u_1^1 = \alpha_1u_1, &
  u_2^1 = \alpha_2u_2 \\
  s_1^1 = \alpha_1s_1, &
  s_2^1 = \alpha_2s_2  \\
  u_1^2 = 2\cdot10^{-7}u_1, &
  u_2^2 = 1.2\cdot10^{-8}u_2, \\
  s_1^2 = 2\cdot10^{-7}s_1, &
  s_2^2 = 2.8\cdot10^{-7}s_2
\end{array}
\end{eqnarray}

 We assume that $H_1^1,H_1^2 \subset \Theta_+$ and
$H_2^1,H_2^2 \subset \Theta_-$. Observe that $I_1\subset H_1^1
\subset H_1^2$ and $I_2\subset H_2^1\subset H_2^2$, where  sets
$I_i$ were defined in the proof of Theorem \ref{thm:exists_per}.
Let
\begin{equation}
W_i=[-\eta_i,\eta_i]\times\{0\}, \quad i=1,2 \label{eq:defW}
\end{equation}
where $\eta_i$ where defined in Theorem \ref{thm:exists_per}. Let
$U_i$, for $i=1,2$ be given by
\begin{equation}
 U_i = H_i^1 + W_i = \{(x+p,\dot{x})\mbox{ } : \mbox{ } (x,\dot{x})\in H_i^1, (p,0)\in W_i\}
 \label{eq:defUi}
\end{equation}

The choices made in (\ref{eq:a}) are motivated by the following
considerations: since we want exploit hyperbolicity of $P$ around
$L_i^*$ it is desirable to choose stable and unstable directions
as $s_i$ and $u_i$. Sets $H_i^1$ (in fact $U_i$) will be used to
establish hyperbolicity around $L_i^*$, hence it desirable to
choose them very small,  as we need to perform a costly rigorous
computation of $DP$ on $U_i$ (a $C^1$-computation). Sets $H_i^2$
are used as a link in the chain of covering relations  between
small $H_i^1$ and relatively large sets $N_i$ defined later in
Section~\ref{sec:homo_hetero}. Since the unstable eigenvalue is
bigger than $10^3$ (see proof of Lemma~\ref{lem:hyperbolic}),  we
can choose $H^2_i$ to be of three order of magnitude larger than
$H^1_i$ and still have a covering relation between them.

The following lemma was proved with a computer assistance
\begin{lem}
\label{lem:dpestm} The maps $P_+:U_1 \to \Theta_+$ and $P_-: U_2
\to \Theta_-$ are well defined. Moreover we have
\begin{equation}
    [DP_+(U_1)]\subset\left(
    \begin{matrix}
        \mathbf A_1 & \mathbf B_1\\
        \mathbf C_1 & \mathbf D_1
    \end{matrix}\right), \qquad
    [DP_-(U_2)]\subset\left(
    \begin{matrix}
        \mathbf A_2 & \mathbf B_2\\
        \mathbf C_2 & \mathbf D_2
    \end{matrix}\right)
\end{equation}
where intervals ${\mathbf A_i}$, ${\mathbf B_i}$, ${\mathbf C_i}$,
${\mathbf D_i}$ are given below
\begin{eqnarray*}
    \begin{array}{cclccl}
    {\mathbf A_1} & = & [695.659,696.1085] &
    {\mathbf B_1} & = & [270.3511,270.4973]\\
    {\mathbf C_1} & = & [1789.9112,1791.46231] &
    {\mathbf D_1} & = & [695.61982,696.12441]\\
    {\mathbf A_2} & = & [573.3983,573.835] &
    {\mathbf B_2} & = & [251.3098,251.4675]\\
    {\mathbf C_2} & = & [1308.1679,1309.5201] &
    {\mathbf D_2} & = & [573.3613,573.848]
    \end{array}
\end{eqnarray*}
\end{lem}

Using the above lemma and symmetry $R$ we can now prove the
following
\begin{lem}
\label{lem:uniq_sym} There exists exactly one fixed point $L_1^*=(x_1^*,0) \in U_1$ for $P_+$.
Moreover we have $|x_1^*-x_1|<\eta_1$.

There exists exactly one fixed point $L_2^*=(x_2^*,0) \in U_2$ for $P_-$. Moreover we have
$|x_2^*-x_2|<\eta_2$.
\end{lem}
{\em Proof:} We write down the proof for $L_1^*$, only. The proof
for $L_2^*$ is analogous.

An easy computation shows that
\begin{displaymath}
  \mbox{det}([DP_+(U_1)] - Id) < 0
\end{displaymath}
hence from Theorem \ref{thm:ewunique} it follows that there exists at most one fixed point for
$P_+$ in $U_1$. Since $I_1 \subset U_1$ then we know from Theorem \ref{thm:exists_per} that one
such fixed point $L_1^*=(x_1^*,0) \in I_1$ exists. The estimate for $|x_1^* - x_1|$ was also given
in Theorem \ref{thm:exists_per}.  \qed

\begin{lem}\label{lem:hyperbolic}
There exist $R$-symmetric h-sets $H_1$ and $H_2$, such that $|H_1|
\subset U_1$, $|H_2| \subset U_2$,t $L_1^* \in H_1$ and $L_2^* \in
H_2$ and the following conditions hold
 \begin{enumerate}
    \item $P_+$ is hyperbolic on $|H_1|$
    \item $P_-$ is hyperbolic on $|H_2|$
 \end{enumerate}
\end{lem}
\noindent
 {\em Proof:} We will proceed as it was outlined in step
2 in section \ref{subsec:howtoprovehomo}. First we need to find a
coordinate frame (via an affine transformation) in which the
inequality (\ref{eq:hyperU}) is satisfied for $(f=P_+,U_f=U_1)$
and $(f=P_-, U_f=U_2)$. Form Lemma~\ref{lem:dpestm} it follows
that $P_+$ is defined on $U_1$ and $P_-$ is defined on $U_2$.

Observe that the transformation of $[DP_+(U_1)]$ ($[DP_-(U_2)]$)
to new coordinates does not depend on the exact location $L_1^*$
($L_2^*)$. In new coordinates $L_1^*=L_2^*=0$, but we have to
choose the coordinate directions in $U_1$ and $U_2$.  It turns out
that the vectors $(u_i,s_i)$ which were used in the definition of
$H^1_i$ are good for this purpose, as they are reasonably good
approximations of unstable  and stable directions of corresponding
Poincar\'e map.

A short  computation shows that in  new coordinates we obtain
\begin{equation}\label{eq:dpestmnewc}
    [DP_+(U_1)] \subset\left(
    \begin{matrix}
        \pmb \lambda_{1,1}  & \pmb \varepsilon_{1,1}\\
        \pmb \varepsilon_{1,2} & \pmb \lambda_{1,2}
    \end{matrix}\right), \qquad
    [DP_-(U_2)]\subset\left(
    \begin{matrix}
        \pmb \lambda_{2,1}  & \pmb \varepsilon_{2,1} \\
        \pmb \varepsilon_{2,2} & \pmb \lambda_{2,2}
    \end{matrix}\right)
\end{equation}
where
\begin{eqnarray*}
    \begin{array}{cclccl}
     \pmb\lambda_{1,1}        & = & [1391.271,1392.239] &
     \pmb\lambda_{1,2}        & = & [-0.482,0.485]\\
     \pmb\varepsilon_{1,1}    & = & [-0.494,0.472] &
     \pmb\varepsilon_{1,2}    & = & [-0.483,0.484]\\
     \pmb\lambda_{2,1}        & = & [1146.751,1147.69] &
     \pmb\lambda_{2,2}        & = & [-0.468,0.47]\\
     \pmb\varepsilon_{2,1}    & = & [-0.481,0.457] &
     \pmb\varepsilon_{2,2}    & = & [-0.468,0.47].
     \end{array}
\end{eqnarray*}

It is clear that $\lambda_{i,2}<1<\lambda_{i,1}$ and
$\varepsilon_{i,1}\varepsilon_{i,2} <
(1-\lambda_{i,2})(\lambda_{i,1}-1)$. Moreover
\begin{eqnarray}\label{eq: cone estimation}
    \frac{\varepsilon_{1,1}}{\lambda_{1,1} - 1}<1<\frac{1-\lambda_{1,2}}{\varepsilon_{1,2}}\\
    \frac{\varepsilon_{2,1}}{\lambda_{2,1} - 1}<1<\frac{1-\lambda_{2,2}}{\varepsilon_{2,2}}
\end{eqnarray}

We define $H_i$ for $i=1,2$ as follows
\begin{equation}
  H_i =t(L_i^*,\alpha_i u_1,\alpha_i s_i),
\end{equation}
where $\alpha_{i}$ were given in (\ref{eq:a}). Observe that by
construction
\begin{equation}
|H_i| \subset |H_i^1| \subset U_i. \label{eq:H-usubU-i}
\end{equation}
 This  shows that $P_+$ is hyperbolic on $|H_1|$ and $P_-$ is hyperbolic on $|H_2|$. \qed

With a computer assistance we proved the following lemma
\begin{lem}
\label{lem:cov_hyper}
 Let $H_1$ and $H_2$ be the h-sets obtained
in Lemma \ref{lem:hyperbolic}, then
\begin{itemize}
  \item $H_1 \cover{P_+} H_1 \cover{P_+} H_1^2$
    \item $H_2^2\cover{P_-} H_2 \cover{P_-} H_2$
\end{itemize}
\end{lem}
{\em Proof:} Consider the following fuzzy h-sets
\begin{equation}
  \widetilde{H_i} = t(W_i,\alpha_i u_i, \alpha_i s_i), \label{eq:defhtilde}
\end{equation}
where $W_i$ are defined by equation (\ref{eq:defW}). We assume
that $|\widetilde{H_1}| \subset \Theta_+$ and $|\widetilde {H_2}|
\subset \Theta_-$.

Observe that $H_i \in \widetilde{H_i}$. The fuzzy sets
$\widetilde{H}_i$ reflect out lack of knowledge of exact
coordinates of $L^*_i$.

\begin{figure}[htbp]
 \centerline{\includegraphics[width=3in]{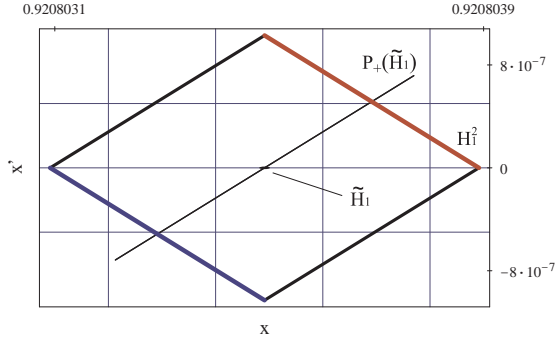}}
 \caption{The Set $H_1^2$ (the large parallelogram), the fuzzy set
${\tilde H}_1$ (a small set in the center) and $P_+(\tilde H_1)$ (the nearly diagonal segment
across $H_1^2$) illustrating covering relation: ${\tilde H}_1\cover{P_+}
\tilde{H}_1\cover{P_+}H_1^2$.
 Vertical edges (when in color:  red and blue) are marked by  a bold line.}
 \label{fig:hyp}
\end{figure}

The following covering relations were established with a computer
assistance (see Fig.~\ref{fig:hyp})
\begin{eqnarray}
     \widetilde{H_1} \cover{P_+} \widetilde{H_1} \cover{P_+}H_1^2 \label{eq:covper1}  \\
    H_2^2 \cover{P_-} \widetilde{H_2} \cover{P_-}
    \widetilde{H_2}  \label{eq:covper2}.
\end{eqnarray}
The assertion of the lemma follows now immediately from Def.
\ref{def:covfuzzy}. \qed

\section{An existence of homo- and heteroclinic connections for Lya\-pu\-nov orbits.}
\label{sec:homo_hetero}
 In this section we prove with a computer assistance
 Theorem~\ref{thm:exists_homo_heteroclinic}. During the proof we
 define h-sets which will be used later in the construction of
 symbolic dynamics in the proof of Theorem~\ref{thm:symb_dyn_1}.

\subsection{An existence of heteroclinic connection between
Lya\-pu\-nov or\-bits.} \label{secsub:hetero} 
\begin{figure}
 \centerline{\includegraphics{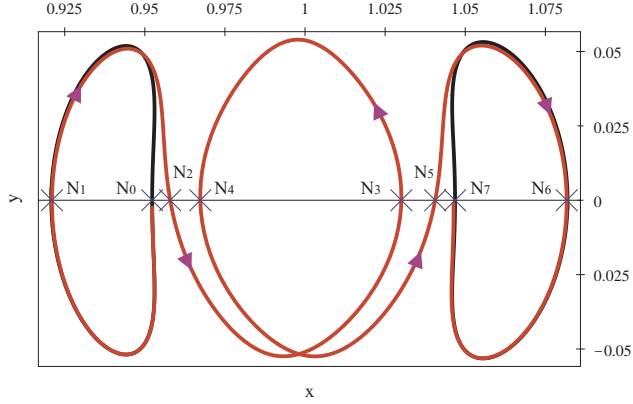}}
 \caption{The location of sets $N_i$ along a  heteroclinic orbit in $(x,y)$-coordinates.}
 \label{fig:covhet}
\end{figure}

In order to prove an
existence of heteroclinic connection between $L_1^*$ and $L_2^*$
we need to find a chain of covering relations which starts close
to $L_1^*$ (begins with $H_1^2$) and ends close to $L_2^*$ (with
$H_2^2$). For this sake we choose the sets $N_i$ along a
numerically constructed, (nonrigorous), heteroclinic orbit in the
vicinity of the intersection of such orbit with the section
$\Theta$ (see Fig. \ref{fig:covhet}).
Let $N_i=t(X_i,u_i,s_i)$ be h-sets, where
\begin{eqnarray*}
 X_0 & = & (0.9522928423486199945,1.23\cdot10^{-5})\\
 X_1 & = & (0.921005737890425169,0.0005205932817646883714)\\
 X_2 & = & (0.957916338594066441,0.02191497366476494527)\\
 X_3 & = & (1.030069865952822683,0.00330658676251664686)\\
 X_4 & = & (0.967306682018305608,0.003703230165036550462)\\
 X_5 & = & (1.040628850444842879,0.02317063455298806404)\\
 X_6 & = & (1.081670357450509545,0.0005918226490172379421)\\
 X_7 & = & (1.046819673646057103,2.13365065043902489\cdot10^{-5})
\end{eqnarray*}
and
\begin{eqnarray*}
\begin{array}{lllclll}
 s_0 & = & (-4\cdot10^{-6},1.45\cdot10^{-5}) & &
 u_0 & = & -R(s_0)/10\\
 s_1 & = & (-4.5\cdot10^{-7},\frac{7}{6}\cdot10^{-6}) & &
 u_1 & = & -R(s_1)/10\\
 s_2 & = & (-1.2\cdot10^{-7},2.92\cdot10^{-7}) & &
 u_2 & = & -R(s_2)\\
 s_3 & = & (-1.05\cdot10^{-7},2.92\cdot10^{-7}) & &
 u_3 & = & -R(s_3)\\
 s_4 & = & (-1\cdot10^{-7},2.9\cdot10^{-7}) & &
 u_4 & = & -R(s_4)/2\\
 s_5 & = & (-1.44\cdot10^{-7},5.8\cdot10^{-7}) & &
 u_5 & = & -R(s_5)/6\\
 s_6 & = & (-1.625\cdot10^{-7},3.75\cdot10^{-7}) & &
 u_6 & = & -R(s_6)/2\\
 s_7 & = & (-8.3\cdot10^{-7},2.9\cdot10^{-6}) & &
 u_7 & = & -R(s_7)/5
\end{array}
\end{eqnarray*}
Vectors $s_i$ were chosen to be a good approximation of the stable
direction at $X_i$. Vectors $u_i$ were chosen as symmetric image
of $s_i$, but  usually with different length.
\begin{rem}
For our computations to succeed vectors $u_i$ might have been
chosen quite arbitrarily (but not too close to $s_i$). For
example, numerical explorations show that if we choose $u_i=\alpha
R(s_i)$ then we can reduce time of rigorous computation
(Lemma~\ref{lem:connect}) around 3 times in comparison to $u_i$
parallel to $x$-axis.
\end{rem}

\begin{figure}[htbp]
  \begin{center}
  \begin{tabular}{cc}
    \includegraphics[width=2in]{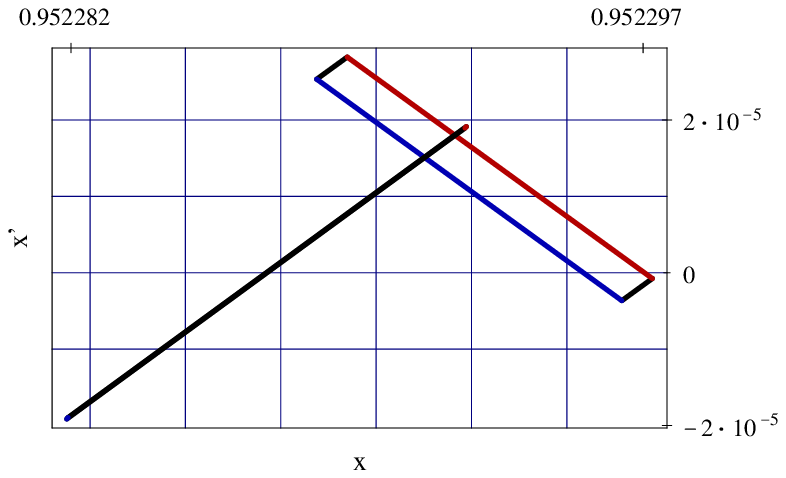} & \includegraphics[width=2in]{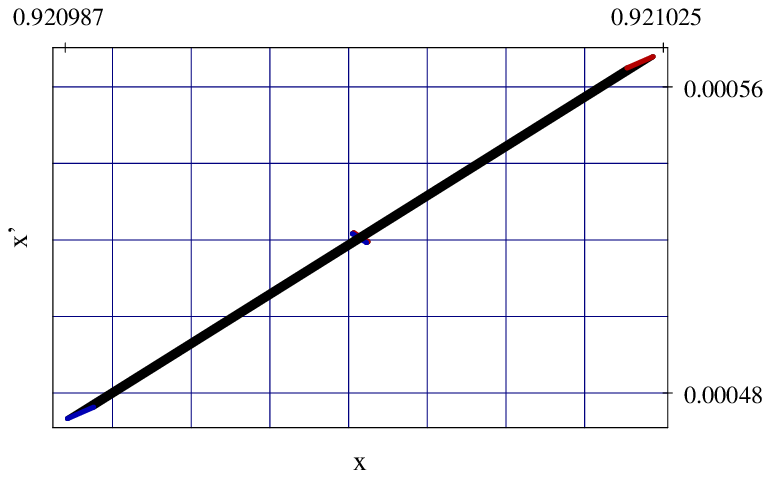} \\
    {\bf\scriptsize a)} & {\bf\scriptsize b)}\\\\
    \includegraphics[width=2in]{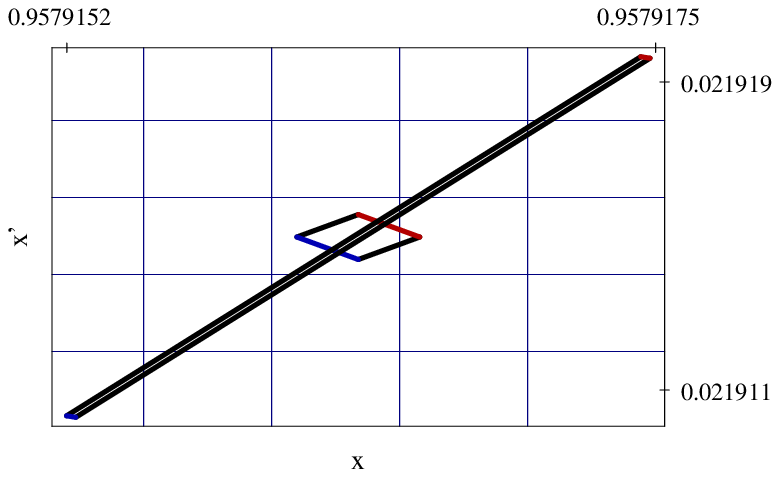} & \includegraphics[width=2in]{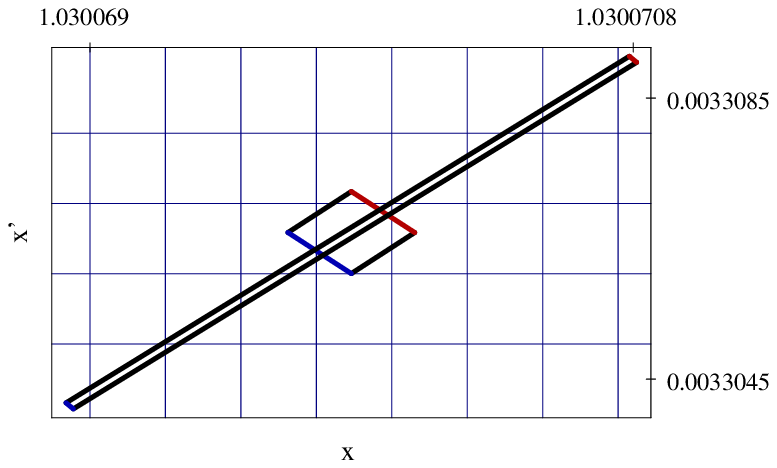} \\
    {\bf\scriptsize c)} & {\bf\scriptsize d)}\\\\
    \includegraphics[width=2in]{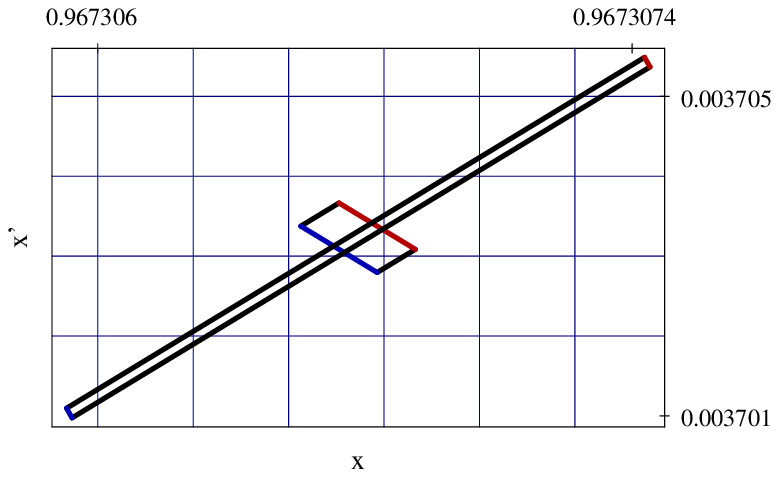} & \includegraphics[width=2in]{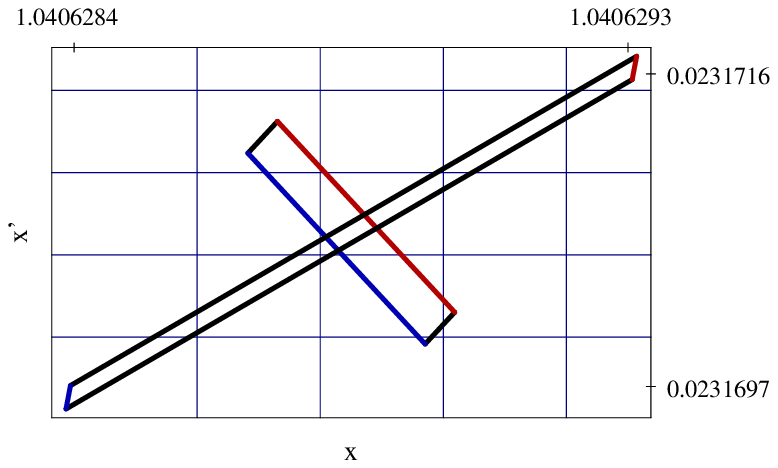} \\
    {\bf\scriptsize e)} & {\bf\scriptsize f)}\\\\
    \includegraphics[width=2in]{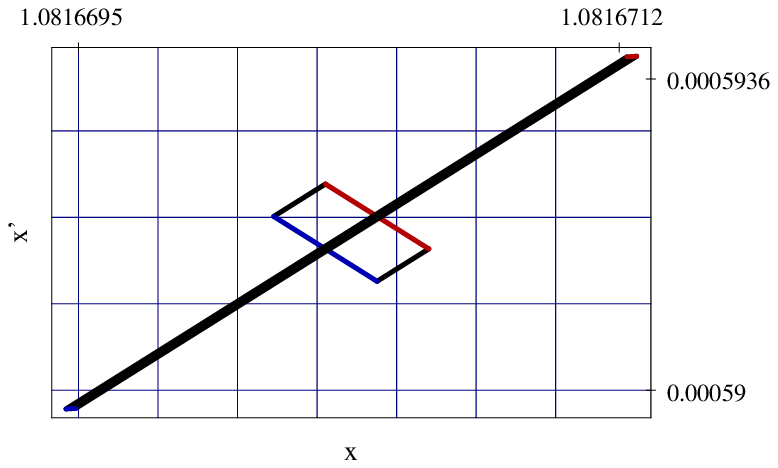} & \includegraphics[width=2in]{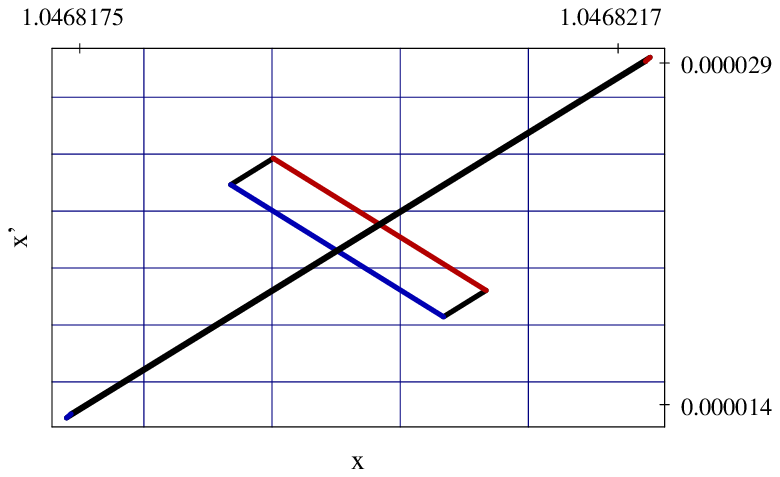} \\
    {\bf\scriptsize g)} & {\bf\scriptsize h)}
 \end{tabular}
 \end{center}
 \caption{A chain of covering relations.
    a) $H_1^2\cover{P_{1/2,+}}N_0$,
    b) $N_0\cover{P_{1/2,-}}N_1$,
    c) $N_1\cover{P_{1/2,+}}N_2$,
    d) $N_2\cover{P_{1/2,-}}N_3$,
    e) $N_3\cover{P_{1/2,+}}N_4$,
    f) $N_4\cover{P_{1/2,-}}N_5$,
    g) $N_5\cover{P_{1/2,+}}N_6$,
    h) $N_6\cover{P_{1/2,-}}N_7$.
    These pictures are'nt produced by a rigorous procedure, as we checked the covering
    relations by  less direct  approach to reduce the computation time -
    see section~\ref{sec:num} for details }
 \label{fig:chain_hetero}
\end{figure}
We assume that $N_0,N_2,N_4,N_6 \subset \Theta_-$ and $N_1,N_3,N_5,N_7 \subset \Theta_+$.
With a~computer assistance we proved
\begin{lem}
\label{lem:connect}
 The maps $$P_{\frac{1}{2},+}:H_1^2\cup N_1\cup N_3\cup N_5\cup N_7\longrightarrow \Theta_-$$ and
 $$P_{\frac{1}{2},-}:N_0 \cup N_2 \cup N_4 \cup N_6\longrightarrow \Theta_+$$  are well defined
 and continuous. Moreover, we have the following chain of covering relations
 \begin{eqnarray*}
  H_1^2\cover{P_{1/2,+}} N_0
  \cover{P_{1/2,-}} N_1
  \cover{P_{1/2,+}} N_2
  \cover{P_{1/2,-}} N_3
  \cover{P_{1/2,+}} N_4\\
  \cover{P_{1/2,-}} N_5
  \cover{P_{1/2,+}} N_6
  \cover{P_{1/2,-}} N_7
  \cover{P_{1/2,+}} H_2^2
 \end{eqnarray*}
\end{lem}
Figure~\ref{fig:chain_hetero} illustrates the chain of covering
relations  from Lemma~\ref{lem:connect} obtained by a  nonrigorous
procedure.

Now we are ready to prove the part of Theorem
\ref{thm:exists_homo_heteroclinic} concerning an existence of
heteroclinic connections.

\begin{thm}
\label{thm:exists_heteroclinic}
 For  PCR3BP with $C=3.03$, $\mu=0.0009537$  there exist two
 periodic solutions in the Jupiter region,$L_1^*$ and $L_2^*$ , called Lyapunov orbits,
 and there exists  heteroclinic connections
 between them, in both directions.
\end{thm}
{\bf Proof:}  We prove only an existence of the connection from
$L_1^*$ to $L_2^*$. An existence of connection in the opposite
direction is obtained by symmetry $R$.

From Lemmas \ref{lem:cov_hyper} and \ref{lem:connect} it follows
that there exists the following chain of covering relations
 \begin{eqnarray*}
  H_1\cover{P_+} H_1
  \cover{P_+}
  H_1^2\cover{P_{1/2,+}} N_0
  \cover{P_{1/2,-}} N_1
  \cover{P_{1/2,+}} N_2
  \cover{P_{1/2,-}} N_3
  \cover{P_{1/2,+}} N_4\\
  \cover{P_{1/2,-}} N_5
  \cover{P_{1/2,+}} N_6
  \cover{P_{1/2,-}} N_7
  \cover{P_{1/2,+}} H_2^2
  \cover{P_-} H_2
  \cover{P_-} H_2
 \end{eqnarray*}
 The assertion follows now from Lemma \ref{lem:hyperbolic} and
 Theorem \ref{th:covhomo}.
\qed

\subsection{Homoclinic connection in an exterior region.}
\label{subsec:hom_ext} 
\begin{figure}[htpd]
   \centerline{\includegraphics[width=2.1in]{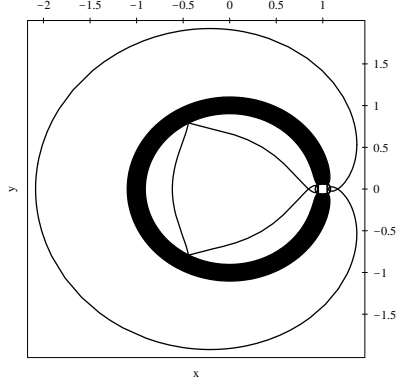}}
   \caption{Homoclinic orbits to $L_1^*$ and $L_2^*$ Lyapunov orbits.\label{fig:homoc}}
\end{figure}

In this section we establish an existence
of an orbit homoclinic to $L_2^*$ (see Figure~\ref{fig:homoc}).
For this end we find a chain of covering relations, which starts
close to $L_2^*$ passes through the sets located in the exterior
region and then ends close to $L_2^*$.  For this sake we choose
the sets $E_i$ along a numerically constructed, (nonrigorous),
homoclinic orbit in the vicinity of the intersection of such orbit
with the section $\Theta$ (see Fig. \ref{fig:homoc}).

We define the following h-sets $E_i=t(Y_i,u_i,s_i)$, where
\begin{eqnarray*}
    Y_0 & = & (-2.08509704964865536,0),\\
    Y_1 & = & (1.160261327316386816,-0.1812035059427922688),\\
    Y_2 & = & (1.059527808809695232,-0.03871458787165545984),\\
    Y_3 & = & (1.082284499686768768,-0.0008090412116073312256),\\
    Y_4 & = & (1.046834433386131072,-0.00002957990840481726976),\\
    Y_5 & = & (1.081929798158888576,-0.0000007068412578518833152)
\end{eqnarray*}
and
\begin{eqnarray*}
\begin{array}{lllclll}
 s_0 & = & (-1\cdot10^{-7},3\cdot10^{-8}) & &
 u_0 & = & -R(s_0)\\
 s_1 & = & (1\cdot10^{-7},8\cdot10^{-8}) & &
 u_1 & = & -4R(s_1)\\
 s_2 & = & (-3\cdot10^{-7},81\cdot10^{-8}) & &
 u_2 & = & -R(s_2)/10\\
 s_3 & = & (-1\cdot10^{-7},23\cdot10^{-8}) & &
 u_3 & = & -R(s_3)/4\\
 s_4 & = & (-1\cdot10^{-7},35\cdot10^{-8}) & &
 u_4 & = & -R(s_4)/4\\
 s_5 & = & (-1\cdot10^{-8},22817915\cdot10^{-15}) & &
 u_5 & = & -R(s_5)/2
\end{array}
\end{eqnarray*}
We assume, that $E_0,E_2,E_4\subset\Theta_+$ and
$E_1,E_3,E_5\subset\Theta_-$.

With a computer assistance we proved the following
\begin{lem}
\label{lem:cov_ext_reg} The maps
\begin{eqnarray*}
  P_{\frac{1}{2},+}&:& E_0 \cup E_2 \cup E_4 \to \Theta_- \\
  P_{\frac{1}{2},-}&:& E_1 \cup E_3 \to \Theta_+ \\
  P_-&:& E_5 \to \Theta_i
\end{eqnarray*}
are well defined and continuous. Moreover, we have the following
chain of covering relations
\begin{eqnarray*}
E_0\cover{P_{1/2,+}}E_1\cover{P_{1/2,-}}E_2\cover{P_{1/2,+}}E_3\cover{P_{1/2,-}}E_4\cover{P_{1/2,+}}E_5\cover{P_-}H_2^2
\end{eqnarray*}
\qed
\end{lem}

We are now ready to state the basic theorem in this section.
\begin{thm}
\label{thm:hom_ext}
 For  PCR3BP with $C=3.03$, $\mu=0.0009537$  there exists a
 an orbit homoclinic to  $L_2^*$.
\end{thm}
{\bf Proof:} From Lemmas~\ref{lem:cov_ext_reg} and
\ref{lem:cov_hyper} it follows that
\begin{equation}
\label{eq:ext_chain}
E_0\cover{P_{1/2,+}}E_1\cover{P_{1/2,-}}E_2\cover{P_{1/2,+}}E_3\cover{P_{1/2,-}}E_4\cover{P_{1/2,+}}E_5\cover{P_-}H_2^2\cover{P_-}H_2.
\end{equation}
Observe that from the definition of $E_0$ it follows that $E_0$ is
$R$-symmetric.  From Corollary~\ref{cor:covsym},
Lemma~\ref{lem:cov_hyper} and equation (\ref{eq:ext_chain}) we
obtain
\begin{eqnarray}
H_2=R(H_2) \invcover{P_-}  R(H_2^2) \invcover{P_-} R(E_5)
\invcover{P_{1/2,-}} R(E_4) \invcover{P_{1/2,+}} \nonumber \\
 R(E_3)\invcover{P_{1/2,-}} R(E_2) \invcover{P_{1/2,+}} R(E_1) \invcover{P_{1/2,-}}
  E_0 = R(E_0) \label{eq:ext_inv}
\end{eqnarray}
From (\ref{eq:ext_chain}), (\ref{eq:ext_inv}),
Lemmas~\ref{lem:cov_hyper} and Theorem~\ref{th:covhomo} we obtain
an orbit homoclinic to $L_2^*$.
 \qed

\subsection{Homoclinic connection in interior region.}
In this section we establish an existence of an orbit homoclinic
to $L_1^*$ (see Figure~\ref{fig:homoc}). For this sake we  find a
chain of covering relations, which starts close to $L_1^*$ passes
through the sets located in the interior region and then ends
close to $L_1^*$. For this sake we choose the sets $F_i$ along a
numerically constructed, (nonrigorous), homoclinic orbit in the
vicinity of the intersection of such orbit with the section
$\Theta$ (see Fig. \ref{fig:homoc}).

We define the following h-sets $F_i=t(Z_i,u_i,s_i)$ where
\begin{eqnarray*}
    Z_0 & = & (-0.6160415155975000064, 0),\\
    Z_1 & = & (0.84668503722876047360, 0.17563753764246766080),\\
    Z_2 & = & (0.94793695784874987520, 0.01522141990729746432),\\
    Z_3 & = & (0.92067611200358768640, 0.00032764933375860776),\\
    Z_4 & = & (0.95228425894935162880, 0.00001048139819208300)
\end{eqnarray*}
and
\begin{eqnarray*}
\begin{array}{lllclll}
 s_0 & = & (-1\cdot10^{-7},25\cdot10^{-8}) & &
 u_0 & = & -R(s_0)\\
 s_1 & = & (1\cdot10^{-7},92\cdot10^{-9}) & &
 u_1 & = & -2.2R(s_1)\\
 s_2 & = & (-25\cdot10^{-9},\frac{33}{4}\cdot10^{-8}) & &
 u_2 & = & -R(s_2)/5\\
 s_3 & = & (-1\cdot10^{-7},26\cdot10^{-8}) & &
 u_3 & = & -R(s_3)/6\\
 s_4 & = & (-1\cdot10^{-7},37\cdot10^{-8}) & &
 u_4 & = & -R(s_4)/6
\end{array}
\end{eqnarray*}
We assume, that $F_0,F_2,F_4\subset\Theta_-$ and
$F_1,F_3\subset\Theta_+$.

With a computer assistance we proved the following
\begin{lem}
\label{lem:cov_inter_reg} The maps
\begin{eqnarray*}
 P_{\frac{1}{2},-}&:& F_0 \cup F_2 \cup F_4 \to \Theta_+ \\
 P_{\frac{1}{2},+}&:& F_1 \cup F_3  \to \Theta_-
\end{eqnarray*}
are well defined and continuous. Moreover, we have the following
covering relations
\begin{eqnarray*}
 F_0\cover{P_{1/2,-}}F_1\cover{P_{1/2,+}}F_2\cover{P_{1/2,-}}F_3\cover{P_{1/2,+}}F_4\cover{P_{1/2,-}}H_1^2
 \end{eqnarray*}
 \qed
\end{lem}
We are now ready to state the basic theorem in this section.
\begin{thm}
\label{thm:hom_inter}
 For  PCR3BP with $C=3.03$, $\mu=0.0009537$  there exists a
 an orbit homoclinic to  $L_1^*$.
\end{thm}
{\bf Proof:} From Lemmas~\ref{lem:cov_inter_reg} and
\ref{lem:cov_hyper} it follows that
\begin{equation}
\label{eq:inter_chain}
 F_0\cover{P_{1/2,-}}F_1\cover{P_{1/2,+}}F_2\cover{P_{1/2,-}}F_3\cover{P_{1/2,+}}F_4\cover{P_{1/2,-}}H_1^2\cover{P_+}H_1
\end{equation}
Observe that from the definition of $F_0$ it follows that $F_0$ is
$R$-symmetric.  From Corollary~\ref{cor:covsym},
Lemma~\ref{lem:cov_hyper} and equation (\ref{eq:inter_chain}) we
obtain
\begin{eqnarray}
 H_1=R(H_1) \invcover{P_+} R(H_1^2) \invcover{P_{1/2,+}} R(F_4)
 \invcover{P_{1/2,-}}R(F_3) \invcover{P_{1/2,+}} R(F_2) \invcover{P_{1/2,-}} \nonumber \\
  R(F_1) \invcover{P_{1/2,+}} R(F_0)=F_0 \label{eq:inter_inv}
\end{eqnarray}
From (\ref{eq:inter_chain}), (\ref{eq:inter_inv}),
Lemmas~\ref{lem:cov_hyper} and Theorem~\ref{th:covhomo} we obtain
an orbit homoclinic to $L_1^*$.
 \qed

\noindent {\bf Proof of
Theorem~\ref{thm:exists_homo_heteroclinic}:} We combine together
Theorems~\ref{thm:exists_heteroclinic}, \ref{thm:hom_ext} and
\ref{thm:hom_inter}.
 \qed

\section{Symbolic dynamics on four symbols }
\label{sec:symbdyn}

The goal of this section is to give a precise meaning  and a proof
to Theorem~\ref{thm:symb_dyn_1}

As in previous sections in the symbol of covering relation we will
drop the degree part, hence we will use $N \cover{f} M$ instead of
$N \cover{f,w} M$ for some nonzero $w$.

From Lemmas~\ref{lem:cov_hyper} and \ref{lem:connect} we know that
there exists the following chain of covering relations
\begin{eqnarray}
  H_1\cover{P_+} H_1
  \cover{P_+}
  H_1^2\cover{P_{1/2,+}} N_0
  \cover{P_{1/2,-}} N_1
  \cover{P_{1/2,+}} N_2
  \cover{P_{1/2,-}} N_3
  \cover{P_{1/2,+}} N_4  \nonumber \\
  \cover{P_{1/2,-}} N_5
  \cover{P_{1/2,+}} N_6
  \cover{P_{1/2,-}} N_7
  \cover{P_{1/2,+}} H_2^2
  \cover{P_-} H_2
  \cover{P_-} H_2. \label{eq:fullchain}
 \end{eqnarray}

From Lemmas~\ref{lem:hyperbolic} and \ref{lem:sym_on_tset} we have
$R(H_i)=H_i$ for $i=1,2$.

From  Lemmas~\ref{lem:cov_hyper}, \ref{lem:connect} and Corollary
\ref{cor:covsym} it follows that
\begin{eqnarray}
  H_2=R(H_2)\invcover{P_-} R(H_2^2)
  \invcover{P_{1/2,-}} R(N_7)
  \invcover{P_{1/2,+}} R(N_6)
  \invcover{P_{1/2,-}} R(N_5) \nonumber \\
  \invcover{P_{1/2,+}} R(N_4)
  \invcover{P_{1/2,-}} R(N_3)
  \invcover{P_{1/2,+}} R(N_2)
  \invcover{P_{1/2,-}} R(N_1) \label{eq:invcovsym} \\
  \invcover{P_{1/2,+}} R(N_0)
  \invcover{P_{1/2,-}} R(H_1^2)
  \invcover{P_{+}} R(H_1)=H_1 \nonumber
\end{eqnarray}

From Lemma~\ref{lem:cov_ext_reg} and the proof of
Theorem~\ref{thm:hom_ext} it follows that
\begin{equation}
\label{eq:ext_chain2}
E_0\cover{P_{1/2,+}}E_1\cover{P_{1/2,-}}E_2\cover{P_{1/2,+}}E_3\cover{P_{1/2,-}}E_4\cover{P_{1/2,+}}E_5\cover{P_-}H_2^2\cover{P_-}H_2.
\end{equation}
and
\begin{eqnarray}
H_2=R(H_2) \invcover{P_-}  R(H_2^2) \invcover{P_-} R(E_5)
\invcover{P_{1/2,-}} R(E_4) \invcover{P_{1/2,+}} \nonumber \\
 R(E_3)\invcover{P_{1/2,-}} R(E_2) \invcover{P_{1/2,+}} R(E_1) \invcover{P_{1/2,-}}
  E_0 = R(E_0). \label{eq:ext_inv2}
\end{eqnarray}

From Lemma~\ref{lem:cov_inter_reg} and the proof of
Theorem~\ref{thm:hom_inter} it follows that
\begin{equation}
\label{eq:inter_chain2}
 F_0\cover{P_{1/2,-}}F_1\cover{P_{1/2,+}}F_2\cover{P_{1/2,-}}F_3\cover{P_{1/2,+}}F_4\cover{P_{1/2,-}}H_1^2\cover{P_+}H_1
\end{equation}
and
\begin{eqnarray}
 H_1=R(H_1) \invcover{P_+} R(H_1^2) \invcover{P_{1/2,+}} R(F_4)
 \invcover{P_{1/2,-}}R(F_3) \invcover{P_{1/2,+}} R(F_2) \invcover{P_{1/2,-}} \nonumber \\
  R(F_1) \invcover{P_{1/2,+}} R(F_0)=F_0. \label{eq:inter_inv2}
\end{eqnarray}

We will  construct now the symbolic dynamics on four symbols. The
construction is a little bit involved, because we have four
different maps in all covering relations listed above.

We assign symbols as follows: $1$ - the set $H_1$, $2$ - $H_2$,
$3$ - $E_0$ and $4$ - $F_0$. The covering relations allow for
transitions $1 \to 1$, $1 \to 2$, $1 \to 4$, $2 \to 1$, $2\to 2$,
$2 \to 3$, $3 \to 2$ and $4 \to 1$. For each such transition $i
\to j$ we associate a pair $(j,i)$. This defines a set of
admissible pairs $\Gamma$.

 For any $(\alpha,\beta) \in \Gamma$ we define a map $f_{(\alpha,\beta)}$ as follows
\begin{displaymath}
f_{(\alpha,\beta)}=  \begin{cases}
    P_+ & \text{if $(\alpha,\beta)=(1,1)$}, \\
  P_- \circ P_{1/2,+}  \circ  (P_{1/2,-} \circ P_{1/2,+})^4 \circ P_+
    & \text{if $(\alpha,\beta)=(2,1)$}, \\
  P_+ \circ P_{1/2,-} \circ  (P_{1/2,+} \circ P_{1/2,-})^4 \circ  P_-
    & \text{if $(\alpha,\beta)=(1,2)$}, \\
  P_-  & \text{if $(\alpha,\beta)=(2,2)$} \\
  P_{\frac{1}{2},+} \circ (P_{1/2,-} \circ P_{1/2,+})^2 \circ
  P_{1/2,+} &  \text{if $(\alpha,\beta)=(4,1)$} \\
  P_{\frac{1}{2},-} \circ (P_{1/2,+} \circ P_{1/2,-})^2 \circ
  P_{-}^2 &  \text{if $(\alpha,\beta)=(3,2)$} \\
  P_+ \circ P_{1/2,-} \circ (P_{1/2,+} \circ P_{1/2,0})^2 &
  \text{if $(\alpha,\beta)=(1,4)$} \\
  P_-^2 \circ P_{1/2,+} \circ (P_{1/2,-} \circ P_{1/2,+})^2 &
    \text{if $(\alpha,\beta)=(2,3)$}
  \end{cases}
\end{displaymath}

Let $\Sigma_\Gamma \subset \{1,2,3,4\}^{\mathbb{Z}}$ be defined as follows $c \in \Sigma_\Gamma$
iff for every $i \in \mathbb{Z}$ $(c_i,c_{i+1}) \in \Gamma$. We set $S_1=H_1$, $S_2=H_2$, $S_3=E_0$
and $S_4=F_0$.

We can now formulate the theorem about an existence of symbolic
dynamics on four symbols symbols
\begin{thm}
\label{thm:symb_dyn_2}
 For any sequence $\alpha=\{\alpha_i\} \in
\Sigma_\Gamma$ there exists a point $x_0 \in S_{\alpha_0}$, such
that
\begin{itemize}
\item its trajectory exists for $t \in (-\infty,\infty)$
 \item $x_n=f_{(\alpha_n,\alpha_{n-1})}\circ \dots \circ
  f_{(\alpha_2,\alpha_1)} \circ f_{(\alpha_1,\alpha_0)} (x_0) \in S_{\alpha_n}$
    for $n >0$
  \item $x_n=f^{-1}_{(\alpha_{n+1},\alpha_{n})}\circ \dots \circ
  f^{-1}_{(\alpha_{-1},\alpha_{-2})} \circ f^{-1}_{(\alpha_{0},\alpha_{-1})} (x_0) \in S_{\alpha_n}$
    for $n < 0$.
\end{itemize}
Moreover, we have
\begin{description}
\item[periodic orbits:] If $\alpha$ is periodic with the principal period equal to
$k$ , then $x_0$ can be chosen so that $x_k=x_0$, hence its
trajectory is periodic.
\item[homo- and heteroclinic orbits:] If $\alpha$ is such that
$\alpha_k=i_-$ for $k \leq k_-$ and $\alpha_k=i_+$ for $k \geq
k_+$, where $i_-,i_+ \in \{1,2\}$, then
\begin{displaymath}
  \lim_{n \to -\infty} x_n = L^*_{i_-}, \qquad \lim_{n \to \infty} x_n = L^*_{i_+}
\end{displaymath}
\end{description}
\end{thm}
{\em Proof:}
 From chains of covering relations
(\ref{eq:fullchain}),(\ref{eq:invcovsym}), (\ref{eq:ext_chain2}),
(\ref{eq:ext_inv2}), (\ref{eq:inter_chain2}),
(\ref{eq:inter_inv2}) and Theorem~\ref{th:top} we obtain the
statement on periodic points for periodic $\alpha$. To treat a
nonperiodic $\alpha$ we approximate it with periodic sequences
$\beta_n$ with increasing periods to obtain sequence of points
$x^n$ and after eventually passing to a subsequence we obtain
$x_0$ with desired properties.

The statement on homo- and heteroclinic orbits is an easy
consequence of Theorem~\ref{th:covhomo} and the hyperbolicity of
$P_\pm$ on $H_i$ established in Lemma~\ref{lem:hyperbolic} \qed

Our methods do not allow to make any claims about the uniqueness
of $x_0$ for a given $\alpha$. The only  claims of this type we
can make is if $\alpha_n=i$ for all $n \in \mathbb{Z}$ then
$x_0=L^*_i$.

\section{Numerical aspects of the proof}
\label{sec:num}

In this section we give details of the computer assisted proofs of
Lemmas \ref{lem:cov_rel_Lap}, \ref{lem:dpestm},
\ref{lem:cov_hyper} and \ref{lem:connect}. As in previous section
in the symbol of covering relation we will drop the degree part,
hence we will use $N \cover{f} M$ instead of $N \cover{f,w} M$ for
some nonzero $w$.

\subsection{The existence and continuity of Poincar\'e maps. Hyperbolicity on $\mathbf U_i$.}

All  proofs required to check first that suitable Poincar\'e maps
($P_\pm$, $P_{\frac{1}{2},\pm}$) are defined on some
parallelograms (supports of our h-sets) on $\Theta_\pm$. For this
end the parallelogram, $Z$, was represented as a finite union of
small parallelograms, $Z_i$, and each of $Z_i$'s was used as an
initial condition for our routine computing the necessary
Poincar\'e map, $P_{1/2,\pm}$ or $P\pm$. We divided horizontal
edges on $n$ equal parts (a horizontal grid) and vertical edges on
$m$ equal parts (a vertical grid) and hence we covered $Z$ by
$n\times m$ parallelograms. Our routine was constructed so that,
if completed successfully, then we can claim that $Z_i$ is
contained in the domain of $P$ and the computed image contains
$P(Z_i)$. Our routine is based on the $C^0$ and $C^1$-Lohner
algorithms \cite{Lo,Z3}.

We had to prove the following assertions
\begin{enumerate}
\item (in Lemma~\ref{lem:cov_rel_Lap}) $P_{\frac{1}{2},+}$ is well defined and continuous on $I_1$ and
  $P_{\frac{1}{2},-}$ is well defined and continuous on $I_2$.
\item (in Lemma~\ref{lem:dpestm}) $P_+$ is well defined and smooth on $U_1$, $P_-$ is well defined
     and smooth on $U_2$.
\item (in Lemma~\ref{lem:cov_hyper} - equations (\ref{eq:covper1}) and (\ref{eq:covper2})).
   $P_-$ is well defined and continuous on $H_2^2$. Observe that since
   $\widetilde H_1\subset U_1$, $\widetilde H_2\subset U_2$,
      then the previous assertion guarantees an existence and
     continuity of $P_+$ on ${\widetilde H}_1$ and $P_-$ on $\widetilde H_2$.
\item (in Lemmas~\ref{lem:connect}, \ref{lem:cov_ext_reg} and \ref{lem:cov_inter_reg})
 $P_{\frac{1}{2},+}$ is well defined and continuous on $H_1^2$, $N_1$, $N_3$,
 $N_5$, $N_7$, $E_0$, $E_2$, $E_4$, $F_1$, $F_3$ . $P_{\frac{1}{2},-}$ is well defined and continuous on
 $N_0$, $N_2$, $N_4$, $N_6$, $E_1$, $E_3$, $F_0$, $F_2$ and $F_4$. $P_-$ is well defined
 on $E_5$.
\end{enumerate}

The first assertion follows easily from the second one. We reason
as follows: since $I_i\subset U_i$, then an existence of $P_-$
($P_+$) on $I_1$ ($I_2$) implies that also $P_{\frac{1}{2},-}$
($P_{\frac{1}{2},+}$) is defined.

To prove the second assertion we cover $U_i$ by finite number ($13\times 13$) of parallelograms.
Then we compute an image of each part and an enclosure of the derivative of the Poincar\'e map
using a routine based on $C^1$-Lohner algorithm recently proposed in \cite{Z3}. As a consequence we
obtain an estimation of $DP_\pm$ (see Lemma~\ref{lem:dpestm}). Parameter settings used in these
computations are listed in Table \ref{tab:hyperbol}. Let us stress also, that a successful
termination of our routine proves also that $P_+$ and $P_-$ are defined on $U_1$ and $U_2$,
respectively. From the standard theory it follows that $P_\pm$ are smooth on their domain.

\renewcommand{\arraystretch}{1.2}
\begin{table}[h]
\begin{center}
 \begin{tabular}{|c|c|c|c|c|}
  \hline
  set & order & step & horizontal grid & vertical grid\\
  \hline
  $U_1$ & 5 & 0.007 & 13 & 13\\
  \hline
  $U_2$ & 5 & 0.007 & 13 & 13\\
  \hline
 \end{tabular}
\end{center}
\caption{Parameter settings of the Taylor method used in $C^1$-computations - in the proof of
Lemma~\ref{lem:dpestm}} \label{tab:hyperbol}
\end{table}

To prove the third and fourth assertion we proceed in the similar way. We cover each set by finite
number of parallelograms and compute an image of each parallelogram. Since an estimation of the
derivative of the Poincar\'e map is not necessary we have used a $C^0$-Lohner algorithm
\cite{Lo,Z3}. Parameter settings for these computations are listed in Table \ref{tab:domain}.

\renewcommand{\arraystretch}{1.6}
\begin{table}[h]
\begin{center}
\begin{tabular}{|c|c|c|c|c|}
   \hline
   covering relations & order & step & h. grid & v. grid \\
   \hline
   \hline
   $\widetilde{H_1}\cover{P_+}H_1^2
   \cover{P_{1/2,+}}N_0
   \cover{P_{1/2,-}}N_1
   \cover{P_{1/2,+}}N_2$ & & & & \\
   $N_2\cover{P_{1/2,-}}N_3
   \cover{P_{1/2,+}}N_4
   \cover{P_{1/2,-}}N_5$ & 5& 0.01&1 &1\\
   $N_5\cover{P_{1/2,+}}N_6
   \cover{P_{1/2,-}}N_7
   \cover{P_{1/2,+}}H_2^2
   \cover{P_-}\widetilde{H_2}$,
     &  &  &  & \\
   \hline
   $E_0\cover{P_{1/2,+}}E_1
   \cover{P_{1/2,-}}E_2
   \cover{P_{1/2,+}}E_3$ & & & &\\
   $E_3\cover{P_{1/2,-}}E_4
   \cover{P_{1/2,+}}E_5
   \cover{P_-}H_2^2$   & 5 & 0.02 & 1 & 1\\
   \hline
   $F_0\cover{P_{1/2,-}}F_1
   \cover{P_{1/2,+}}F_2
   \cover{P_{1/2,-}}F_3$ & & & &\\
   $F_3\cover{P_{1/2,+}}F_4
   \cover{P_{1/2,-}}H_1^2
   \cover{P_+}H_1$   & 5 & 0.02 & 1 & 1\\
   \hline

\end{tabular}
\end{center}
\caption{Parameters of the Taylor method used in the proof of an existence of the Poincar\'e map in
Lemma~\ref{lem:connect}, Lemma~\ref{lem:cov_hyper}, Lemma~\ref{lem:cov_ext_reg} and
Lemma~\ref{lem:cov_inter_reg} \label{tab:domain}}
\end{table}

\subsection{Details of the proof of Lemma~\ref{lem:cov_rel_Lap} }

To prove inequalities (\ref{eq:Lapineq1}),(\ref{eq:Lapineq2}) we
had to compute rigorous enclosures for $P_{\frac{1}{2},+}(x_1 \pm
\eta_1,0)$ and $P_{\frac{1}{2},-}(x_2 \pm \eta_2,0)$,
respectively. The values of the time step and the order of the
Taylor method used in our routine are listed in
Table~\ref{tab:prLap}. Figures \ref{fig:l1} and \ref{fig:l2}
display the actual enclosures obtained.

\renewcommand{\arraystretch}{1.1}
\begin{table}[h]
\begin{center}
 \begin{tabular}{|c|c|c|}
  \hline
  orbit & order of Taylor method & time step\\
  \hline
  $L^*_1$ & 20 & 0.05\\
  \hline
  $L^*_2$ & 19 & 0.055\\
  \hline
 \end{tabular}
\end{center}
\caption{Settings used in the proof of (\ref{eq:Lapineq1}) ($L_1^*$-row) and
(\ref{eq:Lapineq2})($L_2^*$-row) } \label{tab:prLap}
\end{table}

\subsection{How do  we  verify covering relations - details of proofs of
Lemmas \ref{lem:connect}, \ref{lem:cov_ext_reg} and \ref{lem:cov_inter_reg} } This is the most
computationally demanding part of our programm.

In principle the same rigorous computations can be used to obtain
both an existence of Poincar\'e maps and covering relations, but
in practice this does'nt work, i.e. it will result in an enormous
computation time (see the discussion in Sec. 6 of \cite{GZ0}).

It turns out that once an existence of Poincar\'e map is
established, we can reduce the computations to the boundary of our
h-sets and one interval inside, only (see Lemmas~\ref{lem:bd} and
\ref{lem:odwrotne}). Now, when we compute an image of an edge $I$,
we still have to divide it into subintervals, but  the number of
subintervals of the order of square root of the number
parallelograms need to achieve the same accuracy on the
parallelogram build on two intervals of the linear size similar to
that of $I$.

In order to establish an existence of covering relations we need
to verify the assumptions of Theorem~\ref{thm:suffcovu1}.

To facilitate a discussion about various conditions implying
Theorem~\ref{thm:suffcovu1} we introduce the following
\begin{defn}
\label{def:cond} Let $f:\mathbb{R}^{2}\rightarrow\mathbb{R}^{2}$
be a continuous map and let $N_{1}=t(c_1,u_1,s_1)$ and
$N_{2}=t(c_2,u_2,s_2)$ be two h-sets.

 We say that $f$ satisfies condition {\bf ah, a0, a, b+, b-} on $N_1$ and $N_2$ if
\begin{description}
\item[ah:] there exists $q_0 \in [-1,1]$, such that
\begin{displaymath}
f(c_{N_1}([-1,1]\times \{q_0\})) \subset \inte\left(N_{2}^{l}\cup|N_{2}|\cup N_{2}^{r}\right)
\end{displaymath}
\item[a0:] $ f(|N_1|) \cap N^+_2 =\emptyset$
\item[a:]$f(|N_{1}|)\subset\mathrm{int}(N_{2}^{l}\cup|N_{2}|\cup N_{2}^{r})$
\item[b+:] $f(N_{1}^{le})\subset\mathrm{int}(N_{2}^{l})
\quad \mbox{and} \quad
f(N_{1}^{re})\subset\mathrm{int}(N_{2}^{r})$
\item[b-:] $f(N_{1}^{le})\subset\mathrm{int}(N_{2}^{r}) \quad \mbox{and}
\quad f(N_{1}^{re} )\subset\mathrm{int}(N_{2}^{l})$
\end{description}

We say that $f$ satisfies condition {\bf b} on $N_1$ and $N_2$ if
either  {\bf b+} or {\bf b-} is satisfied.
\end{defn}

\begin{rem}
\label{rem:cond} Observe that conditions (\ref{eq:imfu1}),
(\ref{eq:imfu3}), (\ref{eq:covorpr}) and (\ref{eq:covorrev}) from
Theorem~\ref{thm:suffcovu1} coincide with conditions  {\bf ah},
{\bf a0}, {\bf b+} and {\bf b-}, respectively.

Observe that  condition {\bf a} implies conditions {\bf ah} and
{\bf a0}.
\end{rem}

The following lemma gives sufficient conditions for an existence
of covering relations for injective maps.
\begin{lem}
\label{lem:bd} Let $f:\mathbb{R}^{2}\rightarrow\mathbb{R}^{2}$ be
a continuous map and let $N_{1}$ and $N_{2}$ be two h-sets. Assume
that $f$ is  injective  on $|N_{1}|$ and $f$ satisfies condition
{\bf b} on $N_1$ and $N_2$ and the following condition {\bf a'}
\begin{description}
\item [a':]$f(\partial|N_{1}|)\subset\mathrm{int}(N_{2}^{l}\cup
|N_{2}|\cup N_{2}^{r})$.
\end{description}
Then
\begin{displaymath}
  N_1 \cover{f} N_2.
\end{displaymath}
\end{lem}
{\em Proof:} From Remark~\ref{rem:cond} and
Theorem~\ref{thm:suffcovu1} it follows that it is enough to verify
condition  {\bf a}. This follows easily from {\bf a'} and the
Jordan theorem (see \cite{GZ0}, page 180). \qed

Figures illustrating covering relations obtained in Lemmas~\ref{lem:connect} and
\ref{lem:cov_hyper} suggest that condition {\bf a} is satisfied in all relations.
Unfortunately the
verification of condition {\bf a} ( or {\bf a'}) pose the following difficulty:
In the relation
$N_1 \cover{f} N_2$ the set $|N_1|$ is mapped across of $N_2$, without touching its horizontal
edges, but if $|N_2|$ is small then we need a very good estimation of image
of horizontal edges of
$N_1$. This forces us to make a very fine partition of the boundary of $N_1$,
take small time steps
and a high order in the numerical method resulting in very long computation times.

\begin{figure}[htbp]
 \centerline{\includegraphics[width=3in]{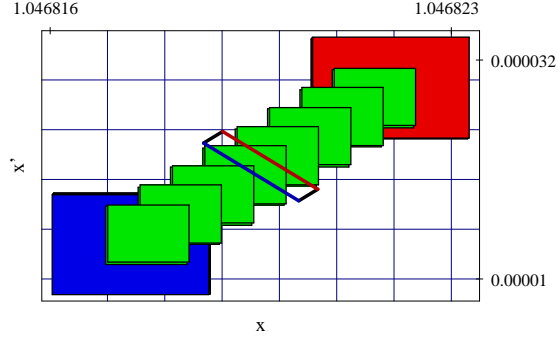}}
 \caption{An example of the rigorous enclosure of an image of $\partial N_6$ in
    relation $N_6\cover{P_{1/2,-}}N_7$.}
 \label{fig:rigorus}
\end{figure}

Above phenomenon is illustrated on Fig.~\ref{fig:rigorus}, which
shows enclosures obtained from our rigorous routines. On this
picture we can see rigorous enclosure for an image of
$P_{1/2,-}(\bd N_6)$. This image was obtained as follows: we
divided the boundary of the set $N_6$ into some number of
subintervals (see Table~\ref{tab:zestawienie}) and computed an
image of each part via $P_{1/2,-}$. This picture shows, that much
tighter enclosure of an image of horizontal edge was required
compared to an enclosure for an image of vertical edges (for
example edge $(N_6)^{le}$ was divided into $8$ equal parts, but
$(N_6)^{be}$ into $5$ equal parts). In  other covering relations
this disproportion was often much  bigger.

To deal with this problem we use the following lemma, in which we
indirectly verify conditions {\bf a0} and {\bf ah} instead of {\bf
a' }. This approach allowed us to  reduce the computation time by
a factor of 5.

\begin{lem}
\label{lem:odwrotne} Let $N_1=t(c_1,u_1,s_1)$,
$N_2=t(c_2,u_2,s_2)$ are h-sets, $f:|N_1|\longrightarrow
\mathbb{R}^2$ an injection of class $\mathcal{C}^1$. Let $\gamma$
be a horizontal line in $|N_1|$ connecting vertical edges given by
\begin{equation}
    \gamma:[-1,1]\ni t \longrightarrow c_1 + t\cdot u_1 \in |N_1|
\end{equation}
Let $g=(g_1,g_2)=c_{N_2}^{-1}\circ f \circ \gamma$. Assume $f$
satisfies condition {\bf b} on $N_1$ and $N_2$ and the following
conditions hold:
\begin{description}
\item [a1] $\frac{dg_1}{dt}(t) \neq 0$ for $t\in(-1,1)$
\item [a2] there exists $t_0\in(-1,1)$ such that $f(\gamma(t_0))\in\mathrm{int}(|N_2|)$
\item [a3] $f^{-1}( N^+_2) \cap |N_1| = \emptyset$
\end{description}
Then $N_1\cover{f}N_2$.
\end{lem}
\noindent {\em Proof:}
 We need to show that conditions {\bf ah} and {\bf a0} are
 satisfied.

Observe that condition {\bf a0} follows immediately from condition
{\bf a3} and injectivity of the map $f$. Namely, by applying $f$
to both sides of {\bf a3} we obtain $ N_2^+ \cap f(|N_1|) =
\emptyset$.

We now show that  condition {\bf ah} is true.

For this end we consider  $f \circ \gamma$ in the coordinates
induced by the map $c_{N_2}$. In these coordinates
\begin{eqnarray}
 |N_2|&=&[-1,1]\times[-1,1]\\
 N_2^r&=&[1,\infty)\times(-\infty,\infty)\\
 N_2^l&=&(-\infty,-1]\times(-\infty,\infty) \\
 f \circ \gamma &=& g
\end{eqnarray}

Without any loss of generality we can assume that
\begin{equation}
\frac{dg_1}{dt}(t)>0, \quad \mathrm{for}\quad t\in(-1,1),
\end{equation}
Hence $g_1$ is a strictly increasing function and from condition {\bf b} it follows that {\bf b+}
is satisfied.

We define two numbers
\begin{eqnarray}
 t^*=\min\{t>t_0\mbox{ }|\mbox{ }f(\gamma(t))\in \partial|N_2|\},\\
 t_*=\max\{t<t_0\mbox{ }|\mbox{ }f(\gamma(t))\in \partial|N_2|\}.
\end{eqnarray}
From conditions {\bf a2}, {\bf a3} and {\bf b} it follows that these numbers are well defined $t_*
< t^*$ and
\begin{eqnarray}
 f(\gamma([-1,t_*))) \subset \mathrm{int}(N_2^l),\\
 f(\gamma((t^*,1]))  \subset \mathrm{int}(N_2^r),\\
 f(\gamma((t_*,t^*)))  \subset \mathrm{int}(|N_2|).
\end{eqnarray}

To finish the proof observe that from condition {\bf b} it follows
that
\begin{displaymath}
  f(\gamma(\pm1)) \in \inte \left(N_2^r \cup  N_2^l\right)
\end{displaymath}
\qed

\begin{rem}
\label{rem:fromdg} Observe that $c_{N_2}^{-1}(x)= A^{-1}(x-c_2)$,
where $A=[u_2^T,s_2^T]$.

Hence
\begin{equation}
\frac{d g_1}{dt}(t)= \sum_{i,j=1}^2 A^{-1}_{1i} df(\gamma(t))_{ij} u_{1,j}, \label{eq:fromdg}
\end{equation}
where $u_1=(u_{1,1},u_{1,2})$.
\end{rem}

\begin{rem}
\label{rem:fromdg1} Observe that, when $N_1=t(c_1,u,s)$  and
$N_2=t(c_2,\alpha u,s_2)$, which means that the unstable
coordinate direction for both h-sets coincide, then
\begin{equation}
  \frac{d g_1}{dt}(t)= \alpha (A^{-1} \cdot df(\gamma(t)) \cdot
  A)_{11}.
\end{equation}
Hence it is enough to look at the $(1,1)$ entry of $df$ expressed
in the coordinates of the h-set $N_2$.
\end{rem}

In Table \ref{tab:zestawienie} we present settings used in the
proof of Lemma~\ref{lem:connect}. In particular the parameter {\em
grid} gives  the number of  equal intervals into which we divide
an edge. A positive time step means that Lemma~\ref{lem:bd} is
used to verify a covering relation. A negative time step means
that we use Lemma~\ref{lem:odwrotne} and symbolizes the fact that
we compute an inverse of the Poincar\'e map to verify condition
{\bf a3}.   Parameter settings for the verification of {\bf a1}
and {\bf a2} are given in Table~\ref{tab:zestawienie1}. In this
table {\em order(m)} and {\em step(m)}  denote an order and a time
step of Taylor method which we use to prove {\bf a1} and  {\em
order(c)}, {\em step(c)} denote an order and a time step of Taylor
method which we use to prove {\bf a2}. Parameter {\em grid(m)}
denotes a number of  equal intervals used to cover  the curve
$\gamma$ in condition {\bf a1}.

To verify condition {\bf a2} we usually compute image of the
center of the set ($t_0=0$). Only in the case
$N_7\cover{P_{1/2,+}}H_2^2$ we used $t_0=0.228$.

\renewcommand{\arraystretch}{1.2}
\begin{table}[htbp]
\begin{center}
   \begin{tabular}{|c|c|c|c|c|}
      \hline
      covering relation & edges & grid & order & step\\
      \hline
      \hline
                           & $(H_1^2)^{be}$ and $(H_1^2)^{te}$ & 6 & 4 & 0.01 \\
      \cline{2-5}
      $H_1^2\cover{P_{1/2,+}}N_0$&$(H_1^2)^{re}$ and $(H_1^2)^{le}$ & 7 & 5 & 0.01\\
      \hline
                           & $N_1^{be}$ and $N_1^{te}$& 1 & 5 & -0.01\\
      \cline{2-5}
      $N_0\cover{P_{1/2,-}}N_1$& $N_0^{re}$ and $N_0^{le}$ & 40 & 8 & 0.04\\
      \hline
                           & $N_2^{be}$ and $N_2^{te}$ & 8 & 6 & -0.01\\
      \cline{2-5}
      $N_1\cover{P_{1/2,+}}N_2$& $N_1^{re}$ and $N_1^{le}$ & 25 & 5 & 0.01\\
      \hline
                           & $N_3^{be}$ and $N_3^{te}$ & 6 & 6 & -0.004\\
      \cline{2-5}
      $N_2\cover{P_{1/2,-}}N_3$& $N_2^{re}$ and $N_2^{le}$ & 5 & 6 & 0.004\\
      \hline
                           & $N_4^{be}$ and $N_4^{te}$ & 3 & 6 & -0.005\\
      \cline{2-5}
      $N_3\cover{P_{1/2,+}}N_4$& $N_3^{re}$ and $N_3^{le}$ & 5 & 6 & 0.004\\
      \hline
                           & $N_5^{be}$ and $N_5^{te}$ & 2 & 6 & -0.01\\
      \cline{2-5}
      $N_4\cover{P_{1/2,-}}N_5$& $N_4^{re}$ and $N_4^{le}$ & 15 & 6 & 0.006\\
      \hline
                           & $N_6^{be}$ and $N_6^{te}$ & 2 & 6 & -0.01\\
      \cline{2-5}
      $N_5\cover{P_{1/2,+}}N_6$& $N_5^{re}$ and $N_5^{le}$ & 32 & 6 & 0.01\\
      \hline
                           & $N_6^{be}$ and $N_6^{te}$ & 8 & 6 & 0.01\\
      \cline{2-5}
      $N_6\cover{P_{1/2,-}}N_7$& $N_6^{re}$ and $N_6^{le}$ & 5 & 6 & 0.01\\
      \hline
                           & $(H_2^2)^{be}$ and $(H_2^2)^{te}$ & 2 & 5 & -0.01\\
      \cline{2-5}
      $N_7\cover{P_{1/2,+}}H_2^2$& $N_7^{re}$ and $N_7^{le}$ & 33 & 5 & 0.01\\
      \hline
   \end{tabular}
\caption{Parameters of the Taylor method used in the proof of Lemma \ref{lem:connect}}
\label{tab:zestawienie}
\end{center}
\end{table}

\renewcommand{\arraystretch}{1.6}
\begin{table}[htbp]
\begin{center}
   \begin{tabular}{|c|c|c|c|c|c|}
      \hline
      covering relations & order(m) & step(m) & grid(m) & order(c) & step(c)\\
      \hline \hline
      $N_0\cover{P_{1/2,-}}N_1\cover{P_{1/2,+}}N_2$     & 4 & 0.01 & 1 & 6 & 0.01\\
      \hline
      $N_2\cover{P_{1/2,-}}N_3$     & 4 & 0.004 & 4 & 6 & 0.004\\
      \hline
      $N_3\cover{P_{1/2,+}}N_4$     & 4 & 0.003 & 2 & 6 & 0.003\\
      \hline
      $N_4\cover{P_{1/2,-}}N_5$     & 4 & 0.004 & 1 & 6 & 0.004\\
      \hline
      $N_5\cover{P_{1/2,+}}N_6\cover{P_{1/2,+}}H_2^2$     & 4 & 0.01 & 1 & 6 & 0.01\\
      \hline
      \hline
      $H_2^2\cover{P_-}\widetilde H_2$   & 4 & 0.01 & 1 &  --- & ---\\
      \hline
   \end{tabular}
\end{center}
\caption{Parameters of the Taylor method used in the proof of
conditions {\bf a1, a2} in Lemmas~\ref{lem:connect} and
\ref{lem:cov_hyper}\label{tab:zestawienie1}}
\end{table}

Parameters of the Taylor method used in the proof of conditions of
Lemma~\ref{lem:cov_ext_reg} are listed in Tables
\ref{tab:cover_exterior}  and \ref{tab:mono_exterior}. In the
relation $E_5\cover{P_-}H_2^2$ we do not compute an image of the
center of set $E_5$ but an image of the point $C=Y_5-0.155u_5$.

\renewcommand{\arraystretch}{1.14}
\begin{table}[htbp]
\begin{center}
   \begin{tabular}{|c|c|c|c|c|}
      \hline
      covering relation & edges & grid & order & step\\
      \hline
      \hline
                                 & $E_1^{be}$ and $E_1^{te}$ & 33 & 7 & -0.04 \\
      \cline{2-5}
      $E_0\cover{P_{1/2,+}}E_1$  & $E_0^{re}$ and $E_0^{le}$ & 12 & 7 & 0.04\\
      \hline
                                 & $E_2^{be}$ and $E_2^{te}$ & 12 & 8 & -0.02\\
      \cline{2-5}
      $E_1\cover{P_{1/2,-}}E_2$  & $E_1^{re}$ and $E_1^{le}$ & 25 & 7 & 0.02\\
      \hline
                                 & $E_3^{be}$ and $E_3^{te}$ & 1 & 7 & -0.02\\
      \cline{2-5}
      $E_2\cover{P_{1/2,+}}E_3$  & $E_2^{re}$ and $E_2^{le}$ & 50 & 7 & 0.03\\
      \hline
                                 & $E_4^{be}$ and $E_4^{te}$ & 2 & 7 & -0.02\\
      \cline{2-5}
      $E_3\cover{P_{1/2,-}}E_4$  & $E_3^{re}$ and $E_3^{le}$ & 7 & 7 & 0.03\\
      \hline
                                 & $E_5^{be}$ and $E_5^{te}$ & 2 & 7 & -0.02 \\
      \cline{2-5}
      $E_4\cover{P_{1/2,+}}E_5$  & $E_4^{re}$ and $E_4^{le}$ & 20 & 7 & 0.03\\
      \hline
                                 & $(H_2^2)^{be}$ and $(H_2^2)^{te}$ & 1 & 5 & -0.02 \\
      \cline{2-5}
      $E_5\cover{P_-}H_2^2$      & $E_5^{re}$ and $E_5^{le}$& 4 & 5 & 0.01\\
      \hline
   \end{tabular}
\caption{Parameters of the Taylor method used in the proof of covering relations in
Lemma~\ref{lem:cov_ext_reg}.\label{tab:cover_exterior}}
\end{center}
\end{table}

\renewcommand{\arraystretch}{1.6}
\begin{table}[htbp]
\begin{center}
   \begin{tabular}{|c|c|c|c|c|c|}
      \hline
      covering relations & o(m) & s(m) & g(m) & o(c) & s(c)\\
      \hline\hline
      $E_0\cover{P_{1/2,+}}E_1\cover{P_{1/2,-}}E_2\cover{P_{1/2,+}}E_3$ & & & & & \\
      $E_3\cover{P_{1/2,-}}E_4\cover{P_{1/2,+}}E_5$     & 4 & 0.02 & 1 & 7 & 0.02\\
      \hline
      $E_5\cover{P_-}H_2^2$         & 5 & 0.02 & 1 & 8 & 0.02\\
      \hline
      \hline
      $F_0\cover{P_{1/2,-}}F_1$     & 5 & 0.01 & 2 & 8 & 0.02\\
      \hline
      $F_1\cover{P_{1/2,+}}F_2\cover{P_{1/2,-}}F_3\cover{P_{1/2,+}}F_4\cover{P_{1/2,-}}H_1^2$     & 4 & 0.01 & 1 & 8 & 0.02\\
      \hline
      $H_1^2\cover{P_+}H_1$         & 4 & 0.01 & 1 & -- & --\\
      \hline
 \end{tabular}
\end{center}
\caption{Parameters of the Taylor method used in the proof of conditions {\bf a1, a2} in
Lemma~\ref{lem:cov_ext_reg} and Lemma~\ref{lem:cov_inter_reg}.
\label{tab:mono_exterior}\label{tab:mono_interior}}
\end{table}

Parameters settings used in the proof of
Lemma~\ref{lem:cov_inter_reg} are listed in Tables
\ref{tab:cover_interior} and \ref{tab:mono_interior}.

\renewcommand{\arraystretch}{1.2}
\begin{table}[htbp]
\begin{center}
   \begin{tabular}{|c|c|c|c|c|}
      \hline
      covering relation & edges & grid & order & step\\
      \hline
      \hline
                                 & $F_1^{be}$ and $F_1^{te}$ & 50 & 6 & -0.01 \\
      \cline{2-5}
      $F_0\cover{P_{1/2,-}}F_1$  & $F_0^{re}$ and $F_0^{le}$ & 20 & 8 & 0.02\\
      \hline
                                 & $F_2^{be}$ and $F_2^{te}$ & 9 & 7 & -0.02\\
      \cline{2-5}
      $F_1\cover{P_{1/2,+}}F_2$  & $F_1^{re}$ and $F_1^{le}$ & 330 & 7 & 0.02\\
      \hline
                                 & $F_3^{be}$ and $F_3^{te}$ & 1 & 7 & -0.02\\
      \cline{2-5}
      $F_2\cover{P_{1/2,-}}F_3$  & $F_2^{re}$ and $F_2^{le}$& 35 & 7 & 0.03\\
      \hline
                                 & $F_4^{be}$ and $F_4^{te}$ & 1 & 7 & -0.02\\
      \cline{2-5}
      $F_3\cover{P_{1/2,+}}F_4$  & $F_3^{re}$ and $F_3^{le}$ & 10 & 7 & 0.03\\
      \hline
                                 & $(H_1^2)^{be}$ and $(H_1^2)^{te}$ & 3 & 7 & -0.02 \\
      \cline{2-5}
      $F_4\cover{P_{1/2,-}}H_1^2$& $F_4^{re}$ and $F_4^{le}$ & 45 & 7 & 0.03\\
      \hline
                                 & $H_1^{be}$ and $H_1^{te}$ & 3 & 8 & -0.02 \\
      \cline{2-5}
      $H_1^2\cover{P_+}H_1$      & $(H_1^2)^{re}$ and $(H_1^2)^{le}$ & 7 & 7 & 0.03\\
      \hline
   \end{tabular}
\caption{Parameters of the Taylor method used in the proof of covering relations for sets
$F_i$.\label{tab:cover_interior}}
\end{center}
\end{table}

\subsection{Verification of covering relations for fuzzy set - details of the proof
 of Lemma~\ref{lem:cov_hyper} }
\label{subsec:covfuzzy}

In this subsection we discuss how we verify covering relations for
fuzzy h-sets. It is convenient to think about a fuzzy h-set
${\widetilde N}$ as an parallelogram with thickened edges. We
define the support, left and right edges and left and right sides
of a fuzzy set ${\widetilde N}$ as follows
\begin{equation*}
\begin{matrix}
  |{\widetilde N}|&=& \langle\bigcup_{M \in {\widetilde N}} |M|\rangle,& &
   \partial {\widetilde N}&=&\bigcup_{M \in {\widetilde N}} \partial |M|,\\
   {\widetilde N}^{le} &=& \langle \bigcup_{M \in {\widetilde N}} M^{le}  \rangle, & &
   {\widetilde N}^{re} &=& \langle \bigcup_{M \in {\widetilde N}} M^{re}  \rangle  \\
  {\widetilde N}^{l} &=& \bigcap_{M \in {\widetilde N}} M^{l}, & &
    {\widetilde N}^{r} &=& \bigcap_{M \in {\widetilde N}} M^{r}.    \\
\end{matrix}
\end{equation*}
where by $\langle Z \rangle$ we denoted a convex hull of the set $Z$. We introduce one more
notation for allowed image of the h-set covering ${\widetilde N}$
\begin{displaymath}
  \mbox{strip}({\widetilde N})= \bigcap_{M \in {\widetilde N}} \inte ( M^l \cup |M|
  \cup M^r)
\end{displaymath}

Lemmas \ref{lem:bd} and \ref{lem:odwrotne} can be easily adopted
to fuzzy h-sets. Namely we have
\begin{lem}
\label{lem:fbd} Let $f:\mathbb{R}^{2}\rightarrow\mathbb{R}^{2}$ be
a continuous map and let ${\widetilde N}_{1}$ and ${\widetilde
N}_{2}$ be two fuzzy h-sets. Assume that $f$ is  injective on
$|N_{1}|$ and the following conditions {\bf af} and {\bf bf} are
satisfied:
\begin{description}
\item[af] $f(|{\widetilde N}_{1}|)\subset  \mathrm{strip}({\widetilde  N}_{2})$

\item[bf] either $f({\widetilde N}_{1}^{le})\subset\mathrm{int}({\widetilde N}_{2}^{l})$ and
$f({\widetilde N}_{1}^{re})\subset\mathrm{int}({\widetilde
N}_{2}^{r})$

or $f({\widetilde N}_{1}^{le})\subset\mathrm{int}({\widetilde N}_{2}^{r})$ and $f({\widetilde N}_{1}^{re}%
)\subset\mathrm{int}({\widetilde N}_{2}^{l})$.
\end{description}
Then ${\widetilde N}_1 \cover{f} {\widetilde N}_2$. \qed
\end{lem}

\begin{lem}
\label{lem:f-odwrotne} Let $\widetilde N_1=t(W_1,u_1,s_1)$, $\widetilde N_2=t(W_2,u_2,s_2)$ be
fuzzy h-sets and $f:|\widetilde N_1|\longrightarrow \mathbb{R}^2$ an injection of class
$\mathcal{C}^1$. Let $\gamma$ be a fuzzy horizontal line in $|\widetilde N_1|$, given by
\begin{equation}
    \gamma:[-1,1]\times W_1 \ni (t,w_1) \longrightarrow w_1 + t\cdot u_1 \in |\widetilde N_1|
\end{equation}
For $w_2 \in W_2$ let $N_{2,w_2}=t(w_2,u_2,s_2)$ and
$$g_{w_2}(t,w_1)=(g_{w_2,1},g_{w_2,2})(t,w_1)=c_{N_{2,w_2}}^{-1} \circ f \circ \gamma(t,w_1).$$

 Assume that $W_1$ is connected,  $f$ satisfies condition {\bf bf} on ${\widetilde N}_1$
and ${\widetilde N}_2$ and the following conditions hold:
\begin{description}
\item [af1] $\frac{dg_{w_2,1}}{dt}(t,w_1) \neq 0$ for $t\in(-1,1)$, $w_1\in W_1$, $w_2\in W_2$
\item [af2] there exists $t_0\in(-1,1)$ and $w_1 \in W_1$ such that $f(\gamma(t_0,w_1))\in\mathrm{int}(|\widetilde N_2|)$
\item [af3] $f^{-1}(\widetilde N^+_2) \cap |\widetilde N_1| = \emptyset$
\end{description}
Then $\widetilde N_1\cover{f}\widetilde N_2$. \qed
\end{lem}

\begin{rem}
\label{rem:ffromdg} Observe that (compare Remark~\ref{rem:fromdg}) that $\frac{d g_{w_2,1}}{dt}$
does not depend on $w_2$ and is given by formula (\ref{eq:fromdg}) with the matrix $A$ depending
only on $u_2$ and $s_2$.
\end{rem}

Let us describe how above conditions {\bf af1, af2, af3} were
verified for relations (\ref{eq:covper1}) and (\ref{eq:covper2}),
which we rewrite below for the convenience of the reader
\begin{eqnarray}
     \widetilde{H_1} \cover{P_+} \widetilde{H_1} \cover{P_+}H_1^2 \label{eq:covper1r}  \\
    H_2^2 \cover{P_-} \widetilde{H_2} \cover{P_-}
    \widetilde{H_2}  \label{eq:covper2r}.
\end{eqnarray}

Let us recall (see Section~\ref{subsec:hyperLyap}) that the h-sets
entering above covering relations are given by
\begin{eqnarray*}
  {\widetilde H}_i&=&t(W_i, \alpha_i u_i, \alpha_i u_i) \\
  H^2_1&=&t(h_1, 2\cdot 10^{-7}u_1, 2\cdot 10^{-7}s_1) \\
  H_2^2&=&t(h_2,1.2\cdot 10^{-8}u_2, 2.8\cdot 10^{-8}s_2)
\end{eqnarray*}

For all covering relations, $N \cover{P_\pm} M$, listed above the
unstable vectors for both h-sets entering the relations are
proportional, hence we can apply  Remark~\ref{rem:fromdg1} and
look on $(1,1)$-entry of $dP_\pm$ expressed in $c_M$-coordinates.

Observe that, since ${\widetilde H}_i \subset U_i$ then from Lemma~\ref{lem:dpestm} and equation
(\ref{eq:dpestmnewc}) we obtain an enclosure for $[DP_\pm(|{\widetilde H}_i)]$ expressed in
${\widetilde H}_i$-coordinates. From inspection of $\lambda_{i,1}$ it follows that condition {\bf
af1} holds for covering relations ${\widetilde H}_i \cover{} {\widetilde H}_i$ and ${\widetilde
H}_1 \cover{} H_1^2$. Since $|H_2^2|$ is not contained in $U_2$ we had to verify condition {\bf
af1} for the relation $H_2^2 \cover{P_-} {\widetilde H}_2$. Parameter settings for this computation
is added as the last row to Table~\ref{tab:zestawienie1}.

Condition {\bf af2} is clearly satisfied with $w_2=L^*_1$ for relations (\ref{eq:covper1}) and
$w_2=L_2^*$ for relations (\ref{eq:covper2}).

Conditions {\bf af3} and {\bf bf} must be verified in direct computations. Parameter settings for
these computations are given in Table \ref{tab:fuzzy}.

It turns out that some inclusions involved in condition {\bf b} can be verified at the same time.
For example, to prove that $\widetilde H_1\cover{P_+}\widetilde H_1\cover{P_+}H_1^2$ we need to
verify condition {\bf b+} for both relations.
\begin{eqnarray}
 P_+(\widetilde H_1^{re})\subset \widetilde H_1^r, & P_+(\widetilde H_1^{le})\subset \widetilde H_1^l,\\
 P_+(\widetilde H_1^{re})\subset (H_1^2)^r, & P_+(\widetilde H_1^{le})\subset (H_1^2)^l.
\end{eqnarray}
It is sufficient to show that
\begin{eqnarray}
 P_+(\widetilde H_1^{re})\subset \widetilde H_1^r \cap (H_1^2)^r\\
 P_+(\widetilde H_1^{le})\subset \widetilde H_1^l \cap (H_1^2)^l
\end{eqnarray}

Similarly, to prove that $H_2^2\cover{P_-}\widetilde
H_2\cover{P_-}\widetilde H_2$ we must verify condition {\bf af3}
for both relations. Namely, we have to check that
\begin{eqnarray}
 P_-^{-1}(\widetilde H^+_2) \cap |\widetilde H_2| = \emptyset,\\
 P_-^{-1}(\widetilde H^+_2) \cap |H_2^2| = \emptyset.\label{eq:disjoint}
\end{eqnarray}
Since $|\widetilde H_2| \subset |H_2^2|$ (compare (\ref{eq:a} and (\ref{eq:defhtilde})) it is
sufficient to show (\ref{eq:disjoint}) only.

\renewcommand{\arraystretch}{1.4}
   \begin{table}[htpd]
   \begin{center}
   \begin{tabular}{|c|c|c|c|c|}
      \hline
      covering relations & edges & grid & order & time step \\
      \hline
      \hline
                           & $(\widetilde H_1)^{be}$ and  $(\widetilde H_1)^{te}$& 2 & 6 & -0.01\\
      \cline{2-5}
      $\widetilde H_1\cover{P_+}\widetilde H_1$ & $(\widetilde H_1)^{re}$ and $(\widetilde H_1)^{le}$ & 3 & 6 & 0.01\\
      \cline{2-5}
      $\widetilde H_1\cover{P_+}H_1^2$ & $(H_1^2)^{be}$ and $(H_1^2)^{te}$ & 2 & 5 & -0.01\\
      \hline
      \hline
                           & $(\widetilde H_2)^{be}$ and $(\widetilde H_2)^{te}$ & 4 & 8 & -0.02\\
      \cline{2-5}
      $H_2^2\cover{P_-}\widetilde H_2$ & $(H_2^2)^{re}$ and $(H_2^2)^{le}$ & 32 & 5 & 0.01\\
      \cline{2-5}
      $\widetilde H_2\cover{P_-}\widetilde H_2$ & $(\widetilde H_2)^{re}$ and $(\widetilde H_2)^{le}$ & 2 & 8 & 0.02\\
      \hline
   \end{tabular}
   \end{center}
   \caption{Parameters of the Taylor method used in the proof of the covering relations for fuzzy sets.
   \label{tab:fuzzy}}
 \end{table}

\section{The {\em pcr3bp} program.}
\label{sec:programm} \comment{ The {\em makefile} is in the {\em
pcr3bp} archive. To link and run our program one need take {\em
pcr3bp} archive and call {\em make}. Then we can run our program -
{\em ./pcr3bp}. Program should be run only in {\bf graphics mode}
Xwindows under Linux (tested on Linux Mandrake 7.1) with
resolution $1024\times 768$ or greater.
 It is possible to run our program under Windows
95/98/NT/2000/XP. In the archive projects one can find the Borland
C++ project. In this case it is necessary to modify the file {\em
setting/setting.h}. One should comment line consists {\bf \#define
LINUX} and uncomment line {\bf \#define WIN95}. Program creates
three output files:
\begin{itemize}
\item {\em report} - some information about verification of covering relations, time of computation.
\item {\em derivative\_1} - DP on $H_1$
\item {\em derivative\_2} - DP on $H_2$
\end{itemize}
}

The program {\bf pcr3bp} is based on {\em CAPD} C++ package \cite{C} developed at the Jagiellonian
University, Krakow, Poland. A whole proof took less 32 minutes on Celeron 1.1Ghz processor 
(compiler Borland C++ 5.02 was used in these computations). The source of our program can be found 
at \cite{W}. Below we briefly describe basic data structures used by our program.

We used the following basic classes:
\begin{description}
\item [class interval] - realizes an interval arithmetic (see \cite{MZ} and references given there). An
interval is a a pair of real numbers, each of the  {\em double}
real type (64 bits).
\item [class intervalVector] - realizes interval arithmetic on
vectors
\item [class intervalMatrix]- realizes interval arithmetic on
matrices
\item [class set] - an abstract class to handle various representations of subsets of $\mathbb R^n$
  used in $C^0$-Lohner algorithm \cite{Lo,MZ,Z3}.
\item [class rec2set] - a class derived from the {\em class set}. This is a class of subsets
   in $\mathbb R^n$ of form $x=x_0+Br+Cr_0$ (doubletons in the terminology introduced in
   \cite{MZ}).
\item [class c1set] - an abstract class  used in
$C^1$-computations via $C^1$-Lohner algorithm \cite{Z3}.
\item [class c1rect] - a class derived from the  {\em class c1set}.  This class represents subsets of
  $\mathbb R^{n}$ as $x=x_0+Br$  and the derivative of the flow with respect to initial conditions
    as $D=D_0+B_0R$.
\item [class c1rect2] - a class derived from the  {\em class c1set}.  This class represents subsets of
  $\mathbb R^{n}$ as $x=x_0+Cr_0+Br$  and the derivative of the flow with respect to initial conditions
    as $D=D_0+B_0R$.
\item [class hSet] - an abstract class representing {\em h-set} in $\mathbb R^n$.
\item [class TripleSet] - a class derived from the  {\em class hSet}. This class implements
  a special case of the {\em hset}, an  h-set on the plane with one stable and unstable
  directions. This class is used to automatize a verification of
  covering relations.
\item [class Taylor] - a class implementing one step of a both $C^0$- and  $C^1$-Lohner algorithms
   for ODEs \cite{Lo,Z3}. It provides the tools which allow to
   compute rigorously the evolution of    sets described by the classes
   given above along the trajectory of an ODE. The Taylor
   coefficients of the solutions of an ODE are generated using an
   {\em automatic differentiation} (see \cite{Lo} and references
   given there).
\item [class PoincareMap] - based on the  {\em class Taylor} provides the tools
to compute the Poincar\'e map and its Jacobian matrix for an ODE
for set classes described above.
\end{description}

The basic methods and operators of the {\em class PoincareMap}
are:
\begin{description}
\item [set\_order(int)] - sets an order of the Taylor method in
the Lohner algorithm.
\item [set\_step(interval)] - sets a time step of the Taylor
method in the Lohner algorithm.
\item [intervalVector operator()(set *)] computes an image of an abstract {\em set} via Poincar\'e map.
\item [intervalVector operator()(c1set *, intervalMatrix \&der)] - computes an image
of an abstract {\em c1set}. On output in {\em der} parameter
contains a matrix of partial derivatives with respect to initial
conditions of the flow, $\frac{\partial \varphi}{\partial x}(t,x)$
for all $x$ defined by parameter {\em c1set} and for  $t \in
[t_1,t_2]$, where $[t_1,t_2]$ contains Poincar\'e return times for
all initial conditions defined by {\em c1set}. From this matrix we
can later  extract easily $\frac{\partial P}{\partial x}$ (see
\cite{Z3} details).
\end{description}

Below we provide a guide to the main functions in our programm.
The {\em main} function is located in the file {\em pcr3bp.cpp}.
Other important  functions are:
\begin{description}
\item [L1\_periodic\_orbit() {\rm and} L2\_periodic\_orbit()] located in the file {\em
lyapunov.cpp}. These functions prove the inequalities in the
Lemma~\ref{lem:cov_rel_Lap}, from which an existence of the
Lyapunov orbits follows.

\item [domain()] located in the file {\em domain.cpp}, used to verify an existence of
the Poincar\'e map for the chain of covering relations form
Lemma~\ref{lem:connect}.

\item [domain\_L1() {\rm and} domain\_L2()] located in the file {\em domain.cpp}, used to
 verify an existence of the Poincar\'e map for covering
relations  $H_1^2\cover{}N_0$ and $H_2^2\cover{}H_2^1$
respectively. ( in Lemmas~\ref{lem:connect} and
\ref{lem:cov_rel_Lap}).

\item [domain\_interior() {\rm and} domain\_exterior()] located in the file {\em domain.cpp},
  used to verify an existence of the Poincar\'e map for covering relations
$F_0\cover{}\ldots\cover{}F_4\cover{}H_1^2\cover{}H_1$ and
$E_0\cover{}\ldots\cover{}E_5\cover{}H_2^2$ (in
Lemmas~\ref{lem:cov_inter_reg} an \ref{lem:cov_ext_reg}).

\item [hyperbolic\_H1() {\rm and} hyperbolic\_H2()] located in file {\em hyperbolic.cpp}. These
functions computes $DP$ on the sets $H_1$ and $H_2$ (see
Lemma~\ref{lem:dpestm})  and proves  covering relations
$\widetilde{H_1}\cover{}H_1^2$ and $H_2^2\cover{}\widetilde{H_2}$
from Lemma~\ref{lem:cov_hyper}.

\item [cover\_heteroclinic()] located in the file {\em cover.cpp}, contains a proof of
Lemma~\ref{lem:connect}, namely it verifies the relations
$H_1^2\cover{}N_0\cover{}\ldots\cover{}N_7\cover{}H_2^2$.

\item [cover\_interior() {\rm and} cover\_exterior()] located in the file {\em
homoclin.cpp}, contains a proof of covering relations from
Lemma~\ref{lem:cov_inter_reg} and Lemma~\ref{lem:cov_ext_reg}.

\item [monotone(), monotone\_interior() {\rm and} monotone\_exterior()]  - located in
the file {\em monotone.cpp}. These functions performs a verification of the monotonicity on the
horizontal central lines connecting vertical edges in \mbox{h-set}. These computations are needed
when we use backward computations. This happens in Lemma~\ref{lem:odwrotne} and
Lemma~\ref{lem:f-odwrotne}.
\end{description}

\section{Concluding remarks, future work}
\label{sec:future} There are several directions in which this
research can be extended.

First, all the methods presented in this paper are not restricted
to the particular parameters of the {\em Oterma} comet, other
parameters may require slight changes in the definition of the
sets on which covering relations are build, but the method itself
will be the same. Basically this method can be applied to prove a
symbolic dynamics in any system for which numerical simulations
indicate an existence of some kind of hyperbolic behavior, for
example here we have homo- and heteroclinic chains.

Another interesting problem is the question of an existence of a
hyperbolic invariant set claimed in \cite{KLMR}, where the authors
assumed an existence of transversal homo- and heteroclinic
connections between Lyapunov orbits and then followed the standard
dynamical system theory argument from the Birkhoff-Smale
homoclinic theorem. Since we did'nt computed here unstable and
stable manifolds, we cannot use these arguments. Observe also that
an rigorous computation of stable and unstable manifolds for our
problem appears to be very difficult  (requires very extensive
$C^1$-computations). Hence developing tools which avoid a direct
computation of invariant manifolds  is of interest. In this
context we formulate the following conjecture.
\begin{con}
Let $f$ be a diffeomorphism. Let $N_0$, $N_1$ be h-sets.
 Assume that $f$ is hyperbolic on $N_0$ and $N_1$ (in the sense
of Def.\ref{def:hyperbolic}). Assume that we have the following
sequences of covering relations
\begin{eqnarray*}
  N_0 \cover{f}  N_0 \cover{f} A_1 \cover{f} A_2 \cover{f} \cdots \cover{f} A_s
  \cover{f} N_1 \\
  N_1 \cover{f}  N_1 \cover{f} B_1 \cover{f} B_2 \cover{f} \cdots \cover{f} B_r
  \cover{f} N_0,
\end{eqnarray*}
then there exists $k \geq 1$ and $S \subset |N_0| \cup |N_1|$,
such that
\begin{itemize}
\item $f^k(S)=S$, i.e. $S$ is an invariant set for $f^k$
\item $S$ is hyperbolic (in the standard sense - see for example \cite{GH})
\item the map $\pi: S \to \Sigma_2=\{0,1\} $ given by $\pi(x)_i =
j$ iff $f^{ki}(x) \in |N_j|$ is one-to-one.
\end{itemize}
\qed
\end{con}

Another interesting problem is a question of stability of obtained
results with respect to various extension  of PCR3BP. By this we
mean the following:
\begin{itemize}
\item Does the symbolic dynamics persists if the Jupiter orbit
become an ellipse with a small eccentricity (which is the case in
nature)? This can be seen as a small periodic perturbation to the
ODE describing PCR3BP. We believe that an answer is positive.
Obviously in this context one can consider a more general
question: \newline
 {\em Assume that we obtained a symbolic
dynamics for an ODE $x'=f(x)$ using covering relations. Does this
symbolic dynamics persists for an nonautonomous ODE $x'=f(x) +
\epsilon(t,x)$ if $\epsilon(t,x)$ is small ? }
\item What about an restricted three body problem in three dimensions? One obvious
observation is that  plane $(x,y)$ is invariant  for full 3D
problem, hence we have symbolic dynamics also in a spatial
problem. We would like to pose a more general question: {\em Does
there exists a symbolic dynamics for 3D problem such the
corresponding orbits are not all contained in the Sun-Jupiter
plane? } Some preliminary numerical explorations in this direction
can be found in paper \cite{GKLMMR},
\item What about full 3-body problem? Will the symbolic dynamics
established here continue to very small but nonzero mass of a
comet?  Some results in this direction for non-degenerate periodic
orbits and generic bifurcations can be found in \cite{MS}.
\end{itemize}

\end{document}